\documentclass{article}
\usepackage[top=1in, bottom=1in, left=1in, right=1in]{geometry}
\usepackage[utf8]{inputenc}
\usepackage{amsmath}
\usepackage{amssymb}
\usepackage{graphicx}
\usepackage{hyperref}
\usepackage{bm}
\usepackage[ruled,vlined]{algorithm2e}
\usepackage{caption}
\usepackage{subcaption}
\usepackage[style=ieee]{biblatex}
\usepackage[english]{babel}
\usepackage{csquotes}
\usepackage[percent]{overpic}
\usepackage{xcolor}
\newtheorem{theorem}{Theorem}[section]

\title{Mobile Sensor Path Planning for Kalman\\ Filter Spatiotemporal Estimation}

\author{Jiazhong Mei$^*$, Steven L. Brunton$^{**}$ and J. Nathan Kutz$^{*,\dag}$\\[.1in] {\small $^*$Department of Applied Mathematics, University of Washington, Seattle, WA 98195}\\ {\small
$^{**}$Department of Mechanical Engineering, University of Washington, Seattle, WA 98105}\\ {\small 
$^\dag$Department of Electrical and Computer Engineering, University of Washington, Seattle WA 98105}\\[.2in]}
\date{}
\begin{document}

\maketitle

\begin{abstract}
    The estimation of spatiotemporal data from limited sensor measurements is a required task across many scientific disciplines. 
    The sensor selection problem, which aims to optimize the placement of sensors, leverages innovations in greedy algorithms and low-rank subspace projection to provide model-free, data-driven estimates. 
    Alternatively, Kalman filter estimation balances model-based information and sparsely observed measurements to collectively make an estimation, with many related optimization algorithms developed for selecting optimal sensors.
    %
    %
    The majority of methods have been developed for stationary sensors, with relatively limited work estimating spatiotemporal data using mobile sensors that leverage both Kalman filtering and low-rank features.
    %
    We show that mobile sensing along dynamic trajectories can achieve the equivalent performance of a larger number of stationary sensors, with performance gains related to three distinct timescales:  (i) the timescale of the spatio-temporal dynamics, (ii) the velocity of the sensors, and (iii) the rate of sampling.
    Taken together, these timescales strongly influence how well-conditioned the estimation task is.  
    Mobile sensing is particularly effective for spatio-temporal data that contain spatially localized structures, whose features  are captured along dynamic trajectories. 
    We draw connections between the Kalman filter performance and the observability of the state space model, and propose a greedy path planning algorithm based on minimizing the condition number of the observability matrix. 
    %
    Through a series of examples of increasing complexity, we show that mobile sensing improves Kalman filter performance in terms of better limiting estimation and faster convergence. 
\end{abstract}


\section{Introduction and Related Work}

Many scientific disciplines require the estimation of spatio-temporal data from limited, point-source sensor measurements for the purpose of characterization, forecasting, reconstructing, and/or controlling a given system. The goal of optimal sensor placement, also known as {\em sensor selection}, is to find the optimal locations in the state space to place only a few sensors so as to achieve the best performance in one or more of the above listed metrics.  
The combinatorial optimization problem of sensor selection is NP-hard, so most algorithms aim to find sub-optimal solutions by leveraging greedy searches and low-rank subspace representations of the system in order to efficiently find a near-optimal solution~\cite{manohar2018data}. Greedy methods~\cite{Tropp:2004,Tropp:2006b} are computationally efficient and include QR decomposition~\cite{trefethen1997numerical} with column pivoting~\cite{manohar2018data, clark2018greedy, saito2021determinant}, (Q)DEIM~\cite{chaturantabut2010nonlinear,sargsyan2015nonlinear,drmac2016new}, and GappyPOD~\cite{everson1995karhunen,astrid2008missing, peherstorfer_stability_2020}, all of which take advantage of the sub-modularity, or near-submodularity, of criteria such as the trace, spectral norm, condition number, determinant and/or its low-rank projection basis. Greedy searches can also be modified to include cost constraints in the sensor placement problem~\cite{clark2018greedy, clark2020cost}. 
Greedy solutions can also be further refined, such as through a genetic algorithm~\cite{sargsyan2018online}.
Other objectives, such as the reconstruction error~\cite{li2021efficient} and the observability matrix~\cite{ilkturk_observability_nodate}, can also be used for sensor selection. Statistical methods using Gaussian process models~\cite{caselton1984optimal,krause2008near} also are effective in leveraging entropy or mutual information as the main objective for optimization. And more recently, shallow decoder networks can be trained within the context of greedy algorithms~\cite{erichson2020shallow,williams2022data}.

In contrast to instantaneous estimation from sensor measurements, Kalman filtering provides a recursive method that estimates based on collective information from prior knowledge of the dynamical model and a time-history of the sensor measurements~\cite{kalman1960,Brunton2019book}. 
In the sensor placement problem, it is often required to have the number of sensors to be at least the same or more than the latent rank of the system in order to be able to capture enough information for reconstruction~\cite{manohar2018data}. 
But with Kalman filter estimation, fewer sensors can be used to achieve the same performance given that the system is observable with these sensor measurements~\cite{Brunton2019book}.  
Commonly, the Kalman filter sensor selection (KFSS) problem studies the objective based on {\it a posteriori} error covariance, which is a metric in Kalman filtering for how much the estimates deviate from the truth. The metric can be considered within an observation period~\cite{tzoumas_sensor_2016}, 
but it is more commonly taken to the limit at the infinite-time horizon when the full convergence of Kalman filtering is reached. Although optimization over the trace of the error covariance matrix, which represents the estimation MSE, does not have a constant-factor polynomial-time approximation~\cite{ye_complexity_2018,zhang_sensor_2017,dhingra_admm_2014},
greedy methods are still near-optimal~\cite{chamon_approximate_2021-1}. 

The diversity of mathematical methods highlighted above for optimal sensor placement typically focus on stationary, point sensors.  However, in many applications, sensors can be mobile,  in which case sensors are allowed to freely move in the measurement space while collecting measurements along the way. The problem concerning the design of trajectories or paths of sensors is called the {\em sensor path planning problem}.
In the field of engineering and robotics, path planning problem has been long considered for the purposes of navigation as well as estimation in a dynamical environment~\cite{gunnarson2021learning, krishna2022finite, biferale2019zermelo, buzzicotti2020optimal, madridano2021trajectory}. 
The task of tracking and estimating a flow field has often been tackled by constructing a simplified, restricted problem that focuses on a network of sensors with a simple formation for efficient paramerization and optimization~\cite{leonard2007collective, devries2013observability, ogren2004cooperative, zhang2010cooperative, paley2020mobile}. 
Different control laws for the path of the sensors are considered for different tasks, including a simple circular or elliptical control~\cite{leonard2007collective}, gradient climbing control~\cite{ogren2004cooperative}, control along level curves~\cite{zhang2010cooperative}, or control based on smoothed particle hydrodynamics~\cite{peng2014dynamic}.
Lynch et al.~\cite{lynch2008decentralized} propose a decentralized mobile network to collectively estimate environmental functions through communication networks, while the sensors move according to a gradient control law that maximizes information. 
Shriwastav et al.~\cite{shriwastav2021dynamic} built a trajectory by connecting a cost-efficient path among optimal sensor placement locations under POD based reconstruction. 
For many of these work, the emphasis is on modeling and control of the sensor positions. 
Sensor scheduling~\cite{liu_optimal_2014,shi_approximate_2013} is a similar problem that concerns a schedule of densely placed sensors. 
Unlike the path planning problem, the sensors do not move in the scheduling problem, although it still can be formulated and solved as a special case of the path planning problem.

While many consider the sensor path in an infinite-time horizon, theoretical studies~\cite{zhang_optimal_nodate,zhao_optimal_2014,mo_infinite-horizon_2014} show that the optimal infinite-time schedule is independent of the initial error covariance and can be approximated arbitrarily closely by a periodic schedule. 
This provides a mathematical foundation for studying problems that consider the planning of a periodic sensor trajectory for spatio-temporal estimation with Kalman filtering. 
Lan and Schwager~\cite{lan_planning_2013,lan_rapidly_2016} approach the periodic path planning problem with a {\em rapidly exploring random cycles} (RRC) method that constructs and evaluates cycles found by randomly exploring the state space using a tree structure;
Chen et al.~\cite{chen_deep_2020} utilize deep reinforcement learning instead as a learnable deterministic method for finding cycles. 
The problem extends to multiple sensors that do not have a set network formation, each following its individual path.  These works are most closely related to the problem considered here.  They approach the combinatorial optimization with a randomized or active search method, first searching for possible cycles, then evaluating their costs. By assumption, the sensors move to a different location at each discrete time step based on the trajectory found in this way. 

\begin{figure}[t!]
    \centering
    \includegraphics[width=0.9\textwidth]{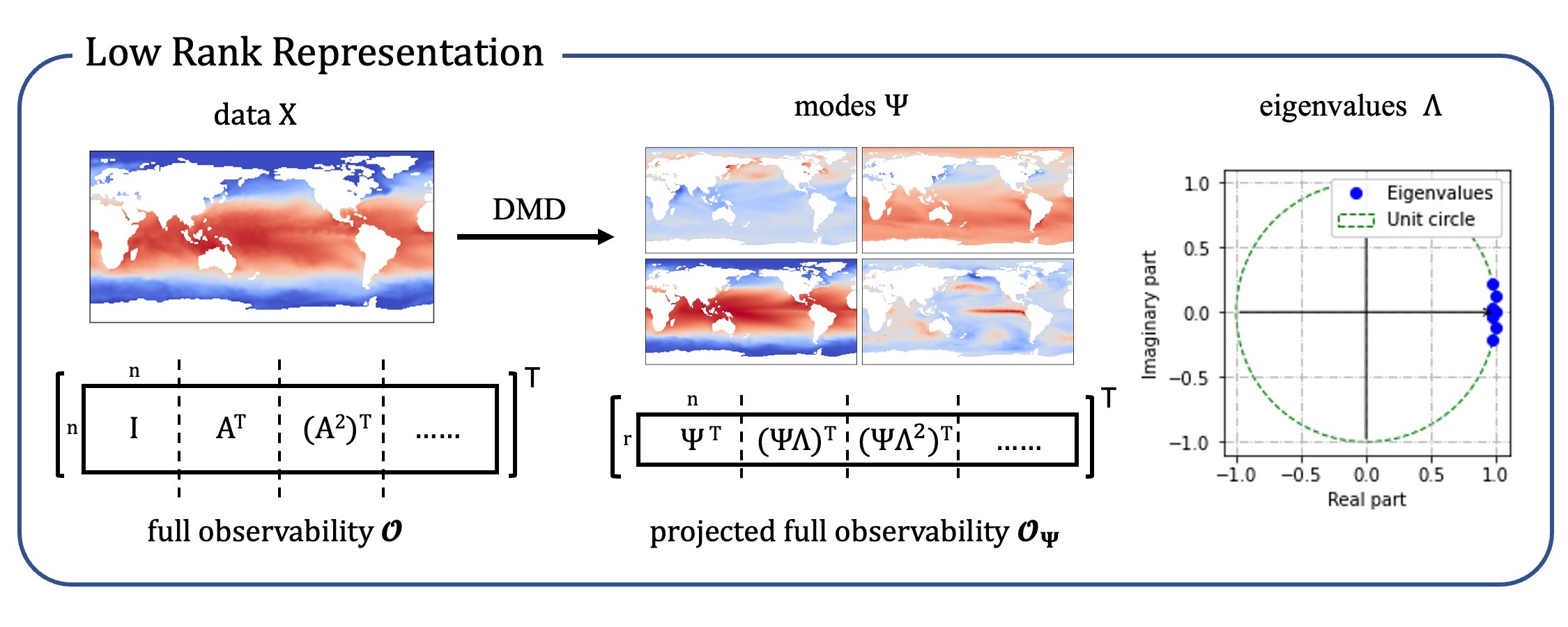} \\
    \includegraphics[width=0.9\textwidth]{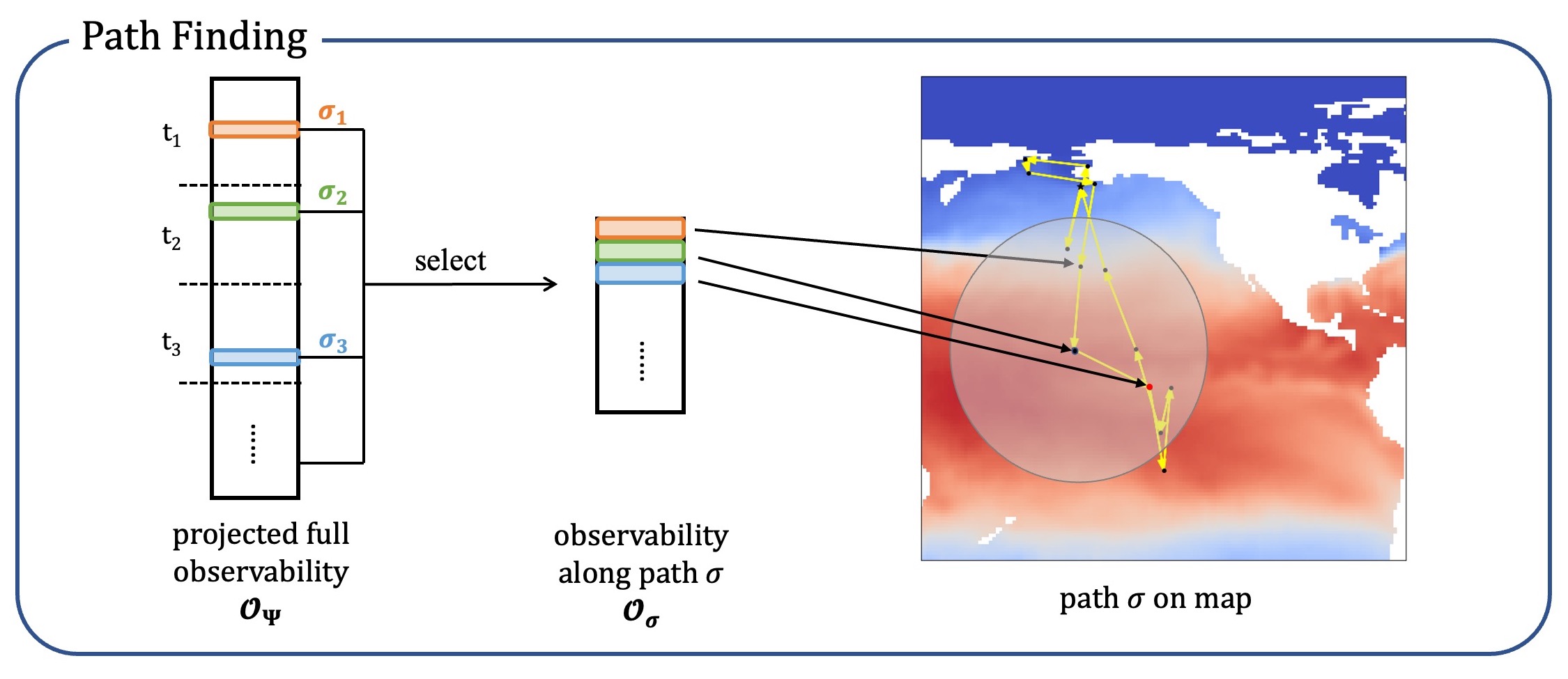}
    \caption{Overview of proposed approach to sensor path planning for dynamic estimation. The panels are divided into two main steps for estimating spatio-temporal data under a Kalman filter setting. 
    The top panel shows the construction of a low-rank representation of the data as the prior model for Kalman filter through DMD. 
    The DMD modes and eigenvalues make up a linear dynamical model in a reduced dimension and a projection back to the original diemsnion. 
    The dimension of the observability matrix is also reduced by the low-rank representation for efficient computation.
    The bottom panel illustrates the greedy path finding algorithm that optimizes the observability matrix along the path and improves Kalman filter estimation performance. 
    It leverages a greedy row selection on the projected full observability matrix.
    Conceptually, at each time step, based on the historical selection of sensor locations, the sensors are led to the next valid locations within a velocity constraint. }
    \label{fig:cover}
\end{figure}

In this paper, we consider the use of mobile sensors to improve the performance of estimating spatio-temporal data with Kalman filtering, where we focus on planning a periodic sensor trajectory that optimizes estimation. 
We consider the condition number of the observability matrix of the model as a metric for the Kalman filter estimation.
The study of observability is not new and has been discussed previously in different sensor problems~\cite{devries2013observability,ilkturk_observability_nodate,manohar_optimal_2021,asghar_complete_2017,rafieisakhaei_use_2017}.
In particular, Manohar et al.~\cite{manohar_optimal_2021} present a balanced model reduction for sensor and actuator selection through observability and controllability in a linear quadratic Gaussian (LQG) controller setting.
We build on these ideas, developing an optimization for the path planning of mobile sensors with the objective of dynamic Kalman filter estimation. 
We identify three distinct timescales related to Kalman filter design and estimation with mobile sensors: 
(i) the timescale of the spatio-temporal dynamics, (ii) the velocity of the sensors, and (iii) the rate of sampling.
We propose an approach for greedy selection based on the empirical observability matrix for path planning, and leverage low-rank representation of the system to promote efficient computation complexity.  Figure~\ref{fig:cover} shows how our overall strategy leverages low-rank approximations in order to determine trajectories in the spatio-temporal fields of interest.
Compared with previous works, our approach does not restrict the formation of the sensor network nor the shape of the trajectory and builds the path by leveraging a low-rank system   representation and greedy optimization. 
Our approach provides a scalable and efficient periodic path planning procedure for multi-sensor and high-dimensional problems.  We conduct a series of experiments on synthetic data, the Kuramoto-Sivashinsky system, and sea surface temperature data to show that mobile sensing improves Kalman filter performance in terms of better limiting estimation and faster convergence. 

\section{Problem Formulation and Background Methods}
The mathematical formulation of the sensor selection problem considers a discrete-time linear system model 
\begin{equation}
    \mathbf{x}_{t+1} = \mathbf{A} \mathbf{x}_t + \mathbf{w}_t,
\end{equation}
where $\mathbf{x}_t \in \mathbb{R}^n$, $\mathbf{A} \in \mathbb{R}^{n \times n}$, and $\mathbf{w}_t \in \mathbb{R}^n$ is the system disturbance following a Gaussian distribution with zero mean and a covariance matrix $\mathbf{0} \prec \mathbf{Q} \in \mathbb{R}^{n \times n}$.
The measurements from $k$ sensors are of the form 
\begin{equation}
    \mathbf{y}_t = \mathbf{C}_t \mathbf{x}_t + \mathbf{v}_t,
\end{equation}
where $\mathbf{y}_t \in \mathbb{R}^k$, $\mathbf{C}_t \in \mathbb{R}^{k \times n}$, and $\mathbf{v}_t \in \mathbb{R}^k$ is the measurement noise following a Gaussian distribution with zero mean and a covariance matrix $\mathbf{0} \prec \mathbf{R} \in \mathbb{R}^{k \times k}$.
Directly measuring in the state space, we write the matrix $\mathbf{C}_t$ as a selection matrix made up of standard unit vectors as columns. 
We further assume that the measurement noise are independent and identical across sensors with variance $\rho$, so the covariance matrix for $\mathbf{v}_t$ is $\mathbf{R} = \rho \mathbf{I}$.
In a time-invariant system with a stationary sensing scenario at fixed locations, then $\mathbf{C}_t = \mathbf{C}$. 

We denote a sensor trajectory $\bm{\sigma} = \{\sigma_1, \sigma_2,...\}$ of $k$ sensors.
$\sigma_t \subseteq [n], |\sigma_t| = k$ is a set containing sensor locations at time $t$, which determine the selection matrix $\mathbf{C}_t = \mathbf{C}(\sigma_t)$ responsible for collecting measurements along the trajectory. 
In general, $\bm{\sigma}$ can extend to an infinite-time horizon. 
Zhang et al. \cite{zhang_optimal_nodate,zhao_optimal_2014} show that any infinite-time trajectory can be approximated by a periodic trajectory. 
Therefore, we focus on a periodic trajectory of fixed cycle rather than a trajectory over infinite-time. 
In practice, periodic trajectories also make sense since many systems contain some periodic or quasi-periodic characteristics. 
Furthermore, it is often favorable to plan a trajectory such that the sensor can return to a specified location periodically for maintenance and sensor recharging.
Then, we write $\bm{\sigma} = \{\sigma_1, \sigma_2,..., \sigma_l\}$ to be a periodic trajectory of length $l$, so that $\sigma_{l+1} = \sigma_1$, $\sigma_{l+2} = \sigma_2$, and so on.

In Sec.~\ref{sec:rom}, we discuss the use of low-rank representation of the system for sparse sampling. 
We then introduce observability of the system  in Sec.~\ref{sec:obs} and relate it to Kalman filter estimation performance in Sec.~\ref{sec:kf}.
Finally, we give attention to three key timescales in the Kalman filter model design in Sec.~\ref{sec:timescale}.

\subsection{Reduced-Order Model and Sparse Sampling} \label{sec:rom}

In order to promote efficient computation and better model representation for sparse sampling, we consider a reduced-order model (ROM).  Specifically, we
consider  a system with a low-rank linear representation,
\begin{equation}
    \begin{split}
    \mathbf{x}_t &= \mathbf{\Psi} \mathbf{z}_t, \quad \mathbf{z}_{t+1} = \mathbf{\Lambda} \mathbf{z}_t + \mathbf{w}_t, \\
    \mathbf{y}_t &= \mathbf{C}_t \mathbf{x}_t + \mathbf{v}_t = \mathbf{C}_t \mathbf{\Psi} \mathbf{z}_t + \mathbf{v}_t,
    \end{split}
    \label{eqn:setup}
\end{equation}
where $\mathbf{z}_t \in \mathbb{R}^m$ ($m < n$) is the internal low-rank dynamics state, $\mathbf{\Lambda} \in \mathbb{R}^{m \times m}$ is the low-rank dynamical system, and $\mathbf{\Psi} \in \mathbb{R}^{n \times m}$ is the linear projection basis.
Measurements $\mathbf{y}_t$ are collected in the original high-dimensional state space. 

One can define a projection basis to be a universal basis for compressed sensing, or a tailored POD basis for a data-driven approach \cite{manohar2018data}. 
However, such basis does not necessarily project to a proper low-rank dynamical system. 
To find a low-rank representation, suppose that the dynamics $\mathbf{A}$ is known, we can take a spectral decomposition of $\mathbf{A}$ and truncate the eigenvalues and eigenvectors to a low-rank representation. 
Alternatively, a data-driven approach is to find a close estimation of the model from the data by using {\em dynamic mode decomposition} (DMD) and its many variants \cite{Tu2014jcd, kutz2016dynamic, jovanovic2014sparsity, askham2018variable,Brunton2022siamreview} that can be useful in sparse sensing \cite{kramer2017sparse}.
DMD modes constitute the linear projection from high-dimensional data to the low-rank representation. 
The DMD eigenvalues form a diagonal dynamics matrix for the low-rank system. 

ROMs are commonly utilized in the stationary sensor placement problem \cite{manohar2018data, drmac2016new,peherstorfer_stability_2020}.
Assuming no disturbance and noise, the measurements can be expressed as 
$\mathbf{y}_t = \mathbf{C}\mathbf{\Psi}\mathbf{z}_t$. Then, we can obtain $\mathbf{x}_t$ through a simple linear reconstruction via the Moore-Penrose pseudoinverse, $\hat{\mathbf{x}}_t = \mathbf{\Psi} \hat{\mathbf{z}}_t = \mathbf{\Psi} (\mathbf{C}\mathbf{\Psi})^\dagger \mathbf{y}_t$.
It is clear that the reconstruction depends on the conditioning of matrix $\mathbf{C}\mathbf{\Psi}$.
Given a tailored basis, Q-DEIM \cite{drmac2016new} uses QR factorization with column pivoting (QRcp) to greedily find near optimal selections.
At each step, QRcp selects a new pivot column with the largest norm and removes the orthogonal projections onto the pivot column from the remaining columns.
Controlling the condition number by maximizing the matrix volume, QRcp enforces a diagonal dominance structure through column pivoting and expands the sub-matrix volume.
In a more recent work, GappyPOD+E \cite{peherstorfer_stability_2020} extends the Q-DEIM method to an ``oversampling'' case where the sample/selection size is larger than the basis rank to improve stability. 
Based on the theory of random sampling in GappyPOD \cite{astrid2008missing}, it is a deterministic method that utilizes a lower bound for the smallest eigenvalue of the submatrix to continue sensor selection over model rank.

\subsection{Observability} \label{sec:obs}

Observability is concerned with the possibility of finding the states of the system from the observations.
A time-varying system of the form $\mathbf{x}_{t+1} = \mathbf{A} \mathbf{x}_t, \mathbf{y}_t = \mathbf{C}_t \mathbf{x}_t$, or a pair $(\mathbf{A}, \mathbf{C}_t)$, is observable at time $t$ if the system state can be determined from the observations in $[t,\tau]$ for some $\tau>t$~\cite{kalman1960general}. 
The system is said to be observable if it is true for all time.
Observability of a system is examined through the observability Gramian.
In our discrete-time system, it is equivalent to study the observability matrix,
$$\bm{\mathcal{O}}_t = \begin{bmatrix} \mathbf{C}_{t} \\ \mathbf{C}_{t+1} \mathbf{A} \\ ... \\ \mathbf{C}_{t+n-1} \mathbf{A}^{n-1} \end{bmatrix}.$$
The system is observable if and only if the observability matrix has full (column) rank. 
When all states are measured, $\mathbf{C}_t = \mathbf{I}$, the full observability matrix is $\bm{\mathcal{O}} = \begin{bmatrix} \mathbf{I} & \mathbf{A}^\top & ... & (\mathbf{A}^{n-1})^\top \end{bmatrix}^\top.$ In the reduced-order model representation, the projected full observability matrix is $$\bm{\mathcal{O}}_\mathbf{\Psi} = \begin{bmatrix} \mathbf{\Psi}^\top & (\mathbf{\Psi} \mathbf{\Lambda})^\top & ... & (\mathbf{\Psi} \mathbf{\Lambda}^{n-1})^\top \end{bmatrix}^\top.$$

In a time invariant system where $\mathbf{C}_t = \mathbf{C}$ is fixed in time, it may need multiple sensors or a long period in time to achieve full rank of the observability matrix. 
For example, for a fully measured system, $\bm{\mathcal{O}}$ is trivially full rank and the system states can be determined immediately at each time step. 
But for sparse sensing on a state space of large dimension, observability of the system is harder to achieve. 
By considering mobile sensing, it opens up possibilities to generate better observability with limited sensors. 

A fully observable system is necessary for an accurate estimation using sparse measurements. 
In particular, it allows Kalman filter estimation to converge to steady-state values on an infinite-time 
horizon. We discuss in more detail the connection between observability and Kalman filter estimation in the following section.

\subsection{Kalman Filter} \label{sec:kf}

We use a Kalman filtering (KF) estimator for spatiotemporal estimation on a linear model. 
Under the assumption of Gaussian noise, KF is known to be the best linear estimator for minimizing mean squared error \cite{kalman1960}. 
Kalman filter algorithms combine the information from the prior knowledge of the system and the observed measurements over time to find an optimal estimate of the system.
Let $\mathbf{\Sigma}_t$ denote the error covariance matrix at time $t$ in the Kalman filter estimation. 
By definition, its trace is the expected squared estimation error at time $t$. 
The error covariance follows a recurrence relation 
$$\mathbf{\Sigma}_{t+1} = \mathbf{A}\mathbf{\Sigma}_t\mathbf{A}^* - \mathbf{A}\mathbf{\Sigma}_t\mathbf{C}_t^*(\mathbf{C}_t\mathbf{\Sigma}_t\mathbf{C}_t^* + \mathbf{R})^{-1}\mathbf{C}_t\mathbf{\Sigma}_t\mathbf{A}^* + \mathbf{Q}. $$
In the time-invariant case, when $t \to \infty$, the limiting error covariance satisfies $\mathbf{\Sigma} = \mathbf{A}\mathbf{\Sigma}\mathbf{A}^* - \mathbf{A}\mathbf{\Sigma}\mathbf{C}^*(\mathbf{C}\mathbf{\Sigma}\mathbf{C}^* + \mathbf{R})^{-1}\mathbf{C}\mathbf{\Sigma}\mathbf{A}^* + \mathbf{Q}$, which is known as the {\em discrete algebraic Riccati equation} (DARE). 
When the system is observable, the error covariance is guaranteed to converge to a limit or a limit cycle in a periodic schedule \cite{zhang_optimal_nodate}.

The limiting error covariance is a common metric to evaluate a KF model. 
However, finding this limit by solving the DARE is difficult and computationally expensive since it does not have a closed-form solution.
Therefore, knowing that observability is a necessary condition for KF estimation performance, we further show how the conditioning of the observability matrix drives the KF estimation.
We first relate the limiting expected squared error from the KF estimation to the conditioning of $\mathbf{C}$ in a time-invariant model with full-rank measurements. We then show that any model with sparse measurements or periodic trajectory can be reformulated at a larger time step to a time-invariant representation with full-rank measurements. And the reformulated selection matrix is the same as the observability matrix in the original form.

The expected squared error is represented as the trace of the error covariance matrix, whose limit is the solution of a DARE. Since the DARE does not have a closed-form solution, we consider an upper and a lower bound for the trace of the solution. While various works have derived different bounds on the DARE solution \cite{dai2011eigenvalue, kwon1996bounds}, we utilize the following results for our analysis:
\begin{theorem}
Consider the DARE $\mathbf{\Sigma} = \mathbf{A}\mathbf{\Sigma} \mathbf{A}^* - \mathbf{A}\mathbf{\Sigma} \mathbf{C}^*(\mathbf{C}\mathbf{\Sigma} \mathbf{C}^* + \mathbf{R})^{-1}\mathbf{C}\mathbf{\Sigma} \mathbf{A}^* + \mathbf{Q}$ with dimension $n$, assuming that $\mathbf{C}^*\mathbf{R}^{-1}\mathbf{C} \succ \mathbf{0}$, $\mathbf{Q} \succ \mathbf{0}$. We then have bounds:

\begin{itemize}
    \item \cite{komaroff1992upper} $tr(\mathbf{\Sigma}) \leq \frac{2tr(\mathbf{Q})}{a_1 + \sqrt{a_1^2 + 4 \lambda_n(\mathbf{C}^*\mathbf{R}^{-1}\mathbf{C}) tr(\mathbf{Q})/n}}$, where $a_1 = 1-\lambda_1(\mathbf{A}^*\mathbf{A}) - \lambda_1(\mathbf{Q}) \lambda_n(\mathbf{C}^*\mathbf{R}^{-1}\mathbf{C})$;
    
    \item \cite{komaroff1992lower} $tr(\Sigma) \geq \frac{2tr(Q^{1/2})^2}{a_2 + \sqrt{a_2^2 + 4n \lambda_1(\mathbf{C}^*\mathbf{R}^{-1}\mathbf{C}) tr(\mathbf{Q}^{1/2})^2}}$, where $a_2 = n-\sum_i |\lambda_i(\mathbf{A})| - tr(\mathbf{Q}^{1/2})^2 \lambda_1(\mathbf{C}^*\mathbf{R}^{-1}\mathbf{C})$.
\end{itemize}

$\lambda_i(\mathbf{X})$ represents the $i$-th largest eigenvalue of $\mathbf{X}$.
\end{theorem}

We can easily show that the lower bound is monotonically decreasing with $\lambda_1(\mathbf{C}^*\mathbf{R}^{-1}\mathbf{C})$, and the upper bound is monotonically decreasing with $\lambda_n(\mathbf{C}^*\mathbf{R}^{-1}\mathbf{C})$ given that $\lambda_1(\mathbf{A}^*\mathbf{A}) \geq 1 - \frac{tr(\mathbf{Q})}{n\lambda_1(\mathbf{Q})}$. 
In the special case $\mathbf{Q}=q\mathbf{I}$, $\lambda_i(\mathbf{Q}) = q$, $tr(\mathbf{Q}) = n\lambda_1(\mathbf{Q})$, $\lambda_1(\mathbf{A}^*\mathbf{A}) \geq 0$ is trivially satisfied. 
The condition is usually satisfied as well for a stable system in general when the disturbance covariance does not have a heavily dominant eigenvalue. 

Since we consider the model with independent and identical measurement noise, $\mathbf{R} = r \mathbf{I}$, so $\lambda_i(\mathbf{C}^*\mathbf{R}^{-1}\mathbf{C})\allowbreak \propto \lambda_i(\mathbf{C}^*\mathbf{C}) = \sqrt{\sigma_i(\mathbf{C})}$, where $\sigma_i(\mathbf{C})$ is the $i$-th largest singular value of $C$. Therefore, in a time-invariant model where $\mathbf{C}$ is full rank, we can minimize the condition number $\kappa(\mathbf{C}) = \frac{\sigma_1(\mathbf{C})}{\sigma_n(\mathbf{C})}$ in order to achieve lower squared error.

However, in most scenarios the system model is more complicated. 
When using limited sensors, the measurements $\mathbf{C}$ will not be full rank.
In the mobile sensor with periodic trajectory scenario where $\mathbf{C}_t$ depends on time, the system is not time-invariant. 
We can show a reformulation of these models to a time-invariant representation in which $\mathbf{C}$ is full rank. Then, the above result applies to these models as well.
Consider the general model $\mathbf{x}_{t+1} = \mathbf{A} \mathbf{x}_t + \mathbf{w}_t$, $\mathbf{y}_t = \mathbf{C}_t \mathbf{x}_t + \mathbf{v}_t$. 
Let $k=n t$ be a larger time step where $n$ is the dimension of the state space or multiples of the sensor trajectory period. Then we can follow \cite{bittanti1991periodic} and rewrite the system as follows:
$$\hat{\mathbf{x}}_{k+1} = \mathbf{x}_{n (t+1)} = \mathbf{A}^n \mathbf{x}_{nt} + \sum_{i=1}^n \mathbf{A}^{i-1}\mathbf{w}_{nt+n-i} := \hat{\mathbf{A}} \hat{\mathbf{x}}_k + \hat{\mathbf{w}}_k,$$
\begin{align*}
\hat{\mathbf{y}}_{k} &:= \begin{bmatrix} \mathbf{y}_{nt} \\ \mathbf{y}_{nt+1} \\ ... \\ \mathbf{y}_{n(t+1)-1} \end{bmatrix} \\
&= \begin{bmatrix} \mathbf{C}_{nt}\mathbf{x}_{nt}+\mathbf{v}_{nt} \\ \mathbf{C}_{nt+1}(\mathbf{A}\mathbf{x}_{nt}+\mathbf{w}_{nt})+\mathbf{v}_{nt+1} \\ ... \\ \mathbf{C}_{n(t+1)-1}(\mathbf{A}^{n-1} \mathbf{x}_{nt} + \sum_{i=2}^{n} \mathbf{A}^{i-2}\mathbf{w}_{nt+n-i})+\mathbf{v}_{n(t+1)-1} \end{bmatrix} \\
&= \begin{bmatrix} \mathbf{C}_{nt} \\ \mathbf{C}_{nt+1}\mathbf{A} \\ ... \\ \mathbf{C}_{n(t+1)-1}\mathbf{A}^{n-1} \end{bmatrix} \mathbf{x}_{nt} + \begin{bmatrix} \mathbf{v}_{nt} \\ \mathbf{C}_{nt+1}\mathbf{w}_{nt}+\mathbf{v}_{nt+1} \\ ... \\ \mathbf{C}_{n(t+1)-1}\sum_{i=2}^{n} \mathbf{A}^{i-2}\mathbf{w}_{nt+n-i}+\mathbf{v}_{n(t+1)-1} \end{bmatrix} \\
&:= \hat{\mathbf{C}} \hat{\mathbf{x}}_k + \hat{\mathbf{v}}_k.
\end{align*}

In the reformation, $\hat{\mathbf{C}}$ is time-invariant, and it is exactly the observability matrix $\bm{\mathcal{O}}$ in the original form. 
If the system is observable, $\bm{\mathcal{O}}$ has full rank, and so does $\hat{\mathbf{C}}$. 
By representing a time-variant system of mobile sensors in a time-invariant form, we can draw the same conclusion as the time-invariant system that the condition number of the observability matrix bounds the limiting error covariance matrix of KF estimation. 
Thus, lowering the condition number of the observability matrix leads to better KF estimation performance.

\subsection{Kalman Filter Design Factors} \label{sec:timescale}
In the mobile sensor scenario, besides planning the trajectory of the sensors, we should also consider in the model design the following key factors: the system Nyquist rate, discrete sampling rate, and sensor speed. 
These three timescales relate to the conditioning of the observability matrix of the system and the performance of Kalman filter estimation. 
Although not a definite guide, the following provides useful heuristics for estimation performance based on these timescales.  

The Nyquist rate represents the internal time scale of the continuous-time dynamics.
It is defined to be twice the highest frequency of the spatiotemporal dynamics. 
The discretization of the continuous-time system is considered good if it samples faster than the Nyquist rate. 
We believe the same applies to mobile sensing with Kalman filter estimation.
At least one measurement should be collected within Nyquist rate to capture the highest frequency feature of the system at the most relevant location. 

The sampling rate refers to the rate at which the measurements are collected. 
It also represents the time step of the discrete model. 
Faster sampling rate above Nyquist adds more measurements in a fixed time interval. In the stationary sensor setting, this improves stability of the estimation. 
With mobile sensors, faster sampling rate further adds more information since the measurement locations change.
This leads to better system observability and KF estimation until the statistical independence of the measurements no longer holds. 

The sensor speed determines the maximum region a sensor can measure in a fixed time interval. 
A faster sensor can reach and observe at a farther location in the state space to achieve better observability. 
More importantly, when the data contains local structures, it is essential to plan the sensor trajectories to capture those structures. 
Faster sensors can achieve this purpose when local structures are well separated in the state space, without the need to assign additional sensors. 

The effect of these timescales will be further discussed in the numerical experiments in Section \ref{sec:num_exp}.

\section{Computing Mobile Sensor Trajectories}

With the problem formulation \ref{eqn:setup}, and the discussion in Section \ref{sec:kf}, we consider the objective to minimize the condition number of the observability matrix under the schedule $\sigma$:
\begin{equation} \label{eqn:obj}
    \min_{\bm{\sigma}: |\bm{\sigma}|=l, |\sigma_i|=k} \kappa(\bm{\mathcal{O}}_{\bm{\sigma}}).
\end{equation}
The observability matrix with respect to trajectory $\bm{\sigma}$ of length $l$ is written as
$$\bm{\mathcal{O}}_{\bm{\sigma}} = \begin{bmatrix} 
\mathbf{C}(\sigma_1) \mathbf{\Psi} \\ 
\mathbf{C}(\sigma_2) \mathbf{\Psi} \mathbf{\Lambda} \\ 
... \\ 
\mathbf{C}(\sigma_l) \mathbf{\Psi} \mathbf{\Lambda}^{l-1} \end{bmatrix}
= \begin{bmatrix} 
\mathbf{C}(\sigma_1) \\ & \mathbf{C}(\sigma_2) \\ && \ddots \\ &&& \mathbf{C}(\sigma_l) \end{bmatrix} \begin{bmatrix} \mathbf{\Psi}  \\ \mathbf{\Psi} \mathbf{\Lambda} \\ ... \\ \mathbf{\Psi} \mathbf{\Lambda}^{l-1} \end{bmatrix} 
:= \mathbf{C}_{\bm{\sigma}} \bm{\mathcal{O}}_\mathbf{\Psi},$$
where $\bm{\mathcal{O}}_\mathbf{\Psi}$ is the projected observability matrix with full measurements, and $\mathbf{C}_{\bm{\sigma}}$ is a block-diagonal selection matrix determined by $\bm{\sigma}$. 
Problem \ref{eqn:obj} becomes a submatrix selection problem minimizing the condition number. 
In the special case when the length of the periodic trajectory is 1, the objective becomes $\max_{\sigma: |\sigma|=k} \kappa(\mathbf{C}_{\bm{\sigma}} \mathbf{\Psi})$, which is identical to that of a stationary sensor placement problem under the DMD basis.
Solving such a problem is in general NP-hard, but just as in the stationary sensor placement problem, we can leverage greedy algorithms and utilize the same idea as QRcp/Q-DEIM for under-sampling and GappyPOD+E or over-sampling to solve it approximately.

We introduce our greedy time-forwarding algorithm in Section \ref{sec:alg} and illustrate on a synthetic example of sparse linear dynamics on a torus in Section \ref{sec:torus} before presenting the main results in Section~\ref{sec:num_exp}.

\subsection{Algorithm} \label{sec:alg}

The projected full observability matrix $\bm{\mathcal{O}}_\mathbf{\Psi}$ is by definition segmented into blocks, so for the purpose of efficient computation, we propose a greedy algorithm that finds sensor locations $\sigma_1, \sigma_2,...,\sigma_l$ by sequentially focusing on each block $\mathbf{\Psi}, \mathbf{\Psi} \mathbf{\Lambda}, ... , \mathbf{\Psi} \mathbf{\Lambda}^{l-1}$.

The algorithms are as follows:

\begin{algorithm}[H]
\SetAlgoLined
\SetKwInput{Input}{Input}
\SetKwInput{Output}{Output}

\Input{low-dimension dynamics matrix $\mathbf{\Lambda}$, projection basis $\mathbf{\Psi}$, number of sensors $k$, trajectory period $l$.}
\Output{trajectory $\bm{\sigma}$.}

Initialize $\mathbf{X} \gets \mathbf{\Psi}$, $\bm{\sigma} \gets []$, $\bm{\mathcal{O}_\sigma} \gets []$\;
\For{$i = 1:l$} {
    
    \For{$j = 1:k$}{
        Find/update candidate rows $S$ in $\mathbf{X}$\;
        Calculate scores and select the next valid row $s \in S$ by Algorithm \ref{alg:score}\;
        Assign the closest and previously unassigned sensor $s' \in \sigma_{i-1}$ to move to $s$, and put $s$ to corresponding position in $\sigma_i$\;
        Append the selected row $\mathbf{X}[s,:]$ to $\bm{\mathcal{O}_\sigma}$\;
    }
    
    Add $\sigma_i$ to $\bm{\sigma}$\;
    Update $\mathbf{X} \gets \mathbf{X\Lambda}$\;
}

\caption{Greedy Time-forwarding Observability-based Path Planning Algorithm}
\label{alg:main}
\end{algorithm}

\begin{algorithm}[H]
\SetAlgoLined
\SetKwInput{Input}{Input}
\SetKwInput{Output}{Output}

\Input{target matrix $\mathbf{X}$, current observability matrix $\bm{\mathcal{O}_\sigma}$, candidates $S$.}
\Output{row index $s$.}


$p \gets$ row size of $\bm{\mathcal{O}_\sigma}$\;
\eIf{$p < m$}{ 
    \tcp{under-sampling, use QRcp rule}
    $[Q,\sim] = qr(\mathcal{O}_\sigma^\intercal)$\;
    $U = Q^\top X^\top$ \;
    $r \gets \sum_{i=p}^m U[i,p:]^2$\;
}{
    \tcp{over-sampling, use GappyPOD+E rule}
    $[\sim,\Sigma,V] = svd(\mathcal{O}_\sigma)$\;
    $g \gets \Sigma_m^2 - \Sigma_{m-1}^2$\;
    $U \gets V^\top X^\top$ \;
    $r \gets g + \sum_{i=1}^m U[i,:]^2$\;
    $r \gets r - \sqrt{r^2 - 4gU[m,:]^2}$\;
}
Sort $r$ in descending order and select $s$ be the first valid ordered index in $S$\;

\caption{Selection Step}
\label{alg:score}
\end{algorithm}

In the selection step, we want to find a row in $\mathbf{X}$ from the candidate set $S$ to append to the current observability matrix in order to minimize $\kappa(\bm{\mathcal{O}_\sigma})$. We use the same selecting rules as in QRcp and GappyPOD+E.
The candidate set $S$ is critical when sensor movement constraints are present. 
When the sensor is unrestricted to move in time, we can simply set $S= [n] \backslash \sigma_i$. 
However, in practice, the sensors have a limit on their speed so there is a restricted region in which the sensors can move between time steps. 
Additionally, the state space can have special multiply-connected topological structure such that not all locations are accessible from every other. 
Future work will incorporate the background flow field in this estimated restriction region, although this is neglected for simplicity in the present work.  

Under a sensor speed constraint, we only consider the locations where 
\begin{itemize}
    \item a sensor can move to within a time step from its current location; 
    \item it can go back to its initial location at the end of the period to form a cycle.
\end{itemize}
These requirements guide the selection of the candidate set $S$ in the algorithm. 
When the topology of the state space is regularly shaped, a simple Euclidean distance can be used; while it is irregular with obstructions or complex network structures, we can resolve to other types of distance functions.

\subsection{Illustrative Example: Sparse Linear Dynamics on a Torus} \label{sec:torus}

To show the effectiveness of mobile sensors, we demonstrate the algorithm with a random simulation of sparse linear dynamical system on a torus~\cite{Brunton2015jcd}. 
We design the system to contain two types of structures: the 2D discrete inverse Fourier transform function and the Gaussian basis function.
A Fourier mode is a global feature present across the state space, while a Gaussian mode is a local feature that only concentrates in a small neighbor around a center. 

On a 128x128 spatial grid, we build the sparse system with 2 conjugate Fourier modes and 3 conjugate Gaussian modes by generating randomly their oscillation frequencies and damping rates. 
This is a system of size $n=128^2=16384$ with a low-dimensional linear representation of rank $m=10$, where the projection basis $\mathbf{\Psi}$ contains the modes, and the low-dimension linear dynamics matrix $\mathbf{\Lambda}$ is diagonal with the oscillation and damping information. 
We generate the data by adding system disturbances and measurement noise.
Since all parameters in the model are known in the synthetic example, we use the trace of the error covariance matrix as an accurate representation of the expected squared error to evaluate the estimation.

\begin{figure}[t]
    \centering
    \hspace{.02\textwidth}
    \begin{overpic}[width=.47\textwidth]{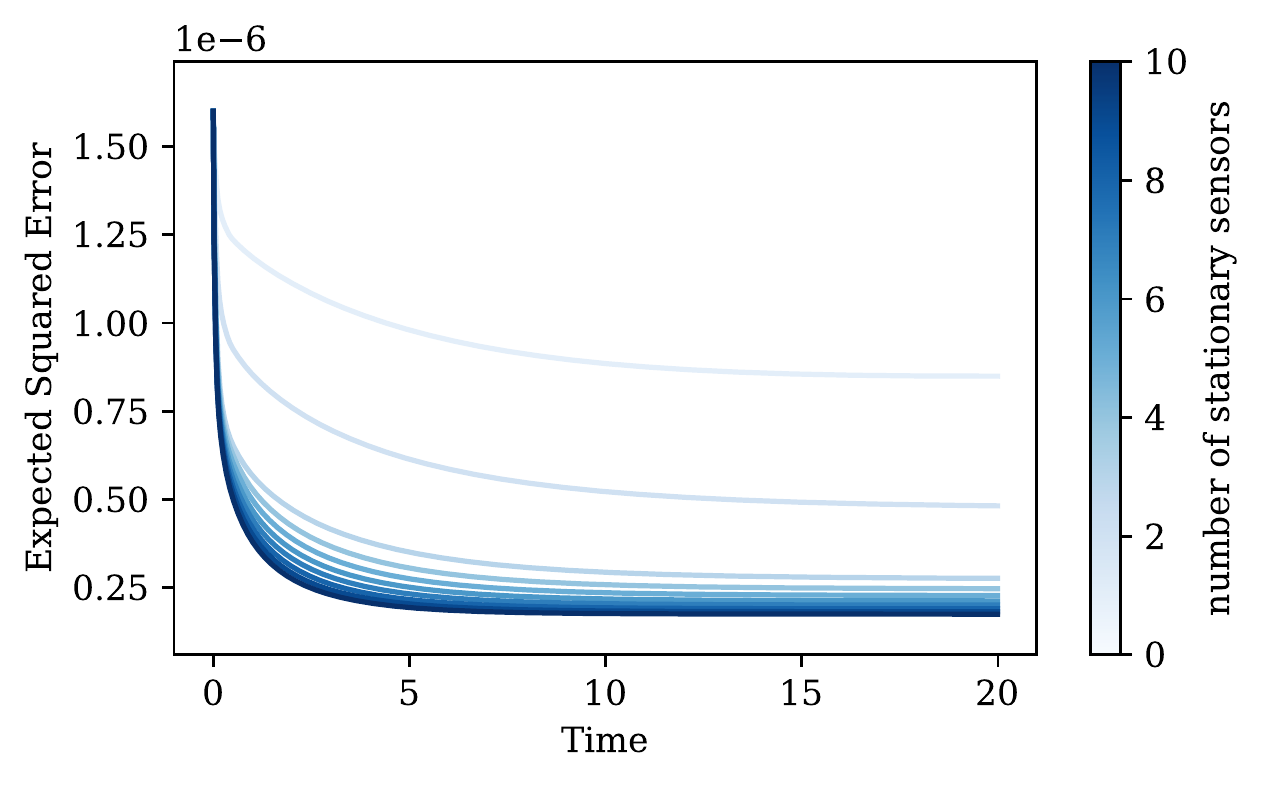}
        \put(-5,60){\color{black}(a)}
    \end{overpic}
    \hspace{.02\textwidth}
    \begin{overpic}[width=.47\textwidth]{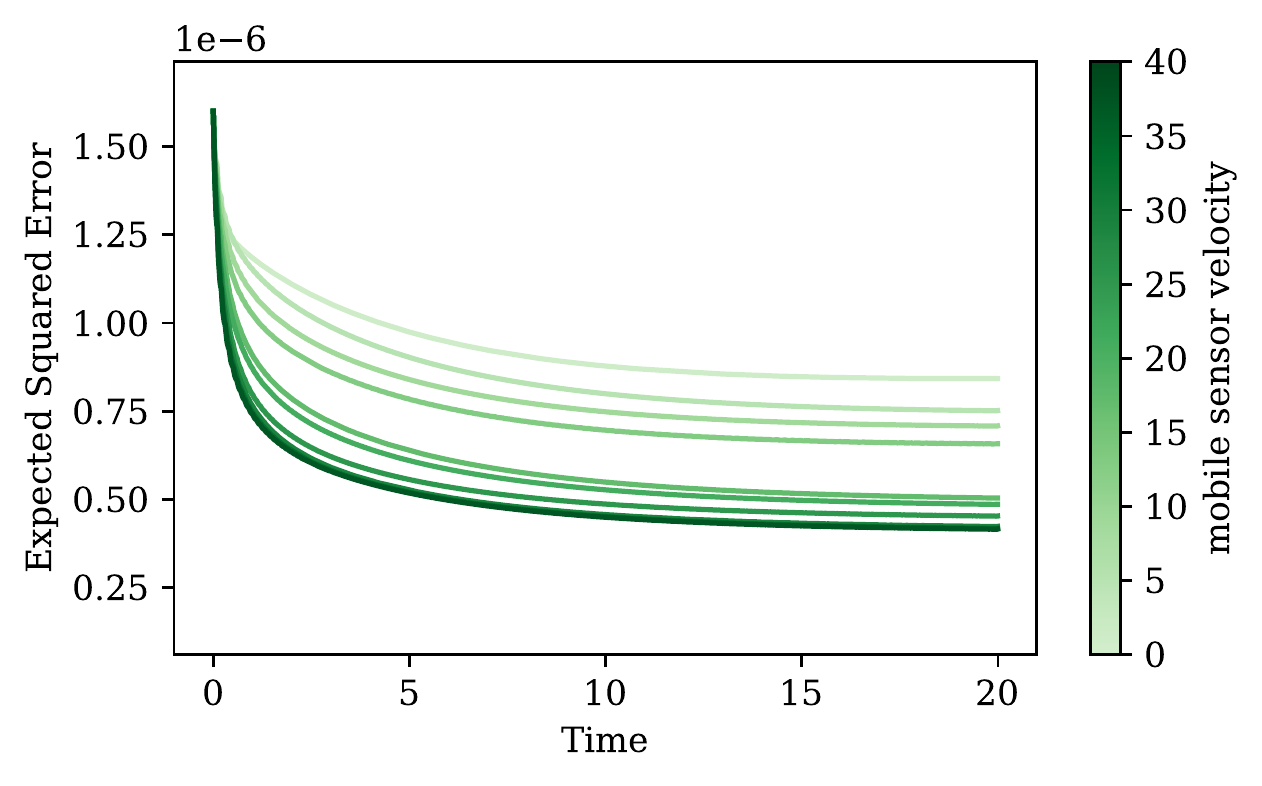}
        \put(-5,60){\color{black}(b)}
    \end{overpic}
     
    \caption{Expected squared error of the KF estimation in time, (a) stationary sensor placement by number of sensors; (b) one mobile sensor by velocity constraints.}
    \label{fig:torus_results}
\end{figure}

First, we estimate the system with sensors at fixed locations selected by applying QRcp on the basis $\mathbf{\Psi}$, a common sensor placement strategy. 
We see from Figure \ref{fig:torus_results} there is significant performance improvement as we increase the number of fixed sensors up to 3. 
At least 3 sensors are needed to obtain a good estimation of the system so that they can be placed to observe the local regions of the Gaussian modes.

\begin{figure}[t]
    \centering
    \includegraphics[height=2.1in]{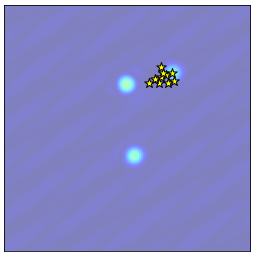}
    \hspace{.1in}
    \includegraphics[height=2.1in]{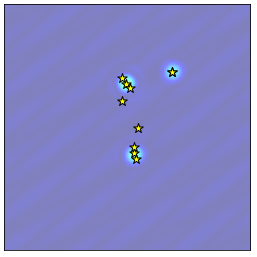}
    
    \caption{Sensor locations along trajectory with speed of 5 (left) and 37 (right) units per time step of 0.01.}
    \label{fig:torus_locs}
\end{figure}

We then show that equivalent performance can be achieve using only one mobile sensor with the same sampling rate and fast enough speed. 
We choose a trajectory period such that the cycle is complete within the Nyquist rate of the system.
When the sensor is slow, there is no significant improvement since the sensor cannot move to other local features within a cycle. 
But, with fast enough speed, our algorithm is able to direct the sensor to reach the localization of all three Gaussian modes and make better estimation (Figure \ref{fig:torus_locs}). 
Under the same sampling rate, three sensors collect three times as many measurements as only one sensor within any time interval. 
This fundamental differences in measurement size due to  number of sensors contributes to the difference between three stationary sensors and one mobile sensors. 
We can narrow this performance gap by increasing the sampling rate. 
At a sampling rate of 0.001, the difference in estimation error is minimal.

Through this synthetic experiment we see that a mobile sensor can indeed improve Kalman filter estimation, and the trajectory planned by our greedy method is effective to pinpoint local structures and achieve good observability.

\begin{figure}[t]
    \centering
    \includegraphics[height=2.1in]{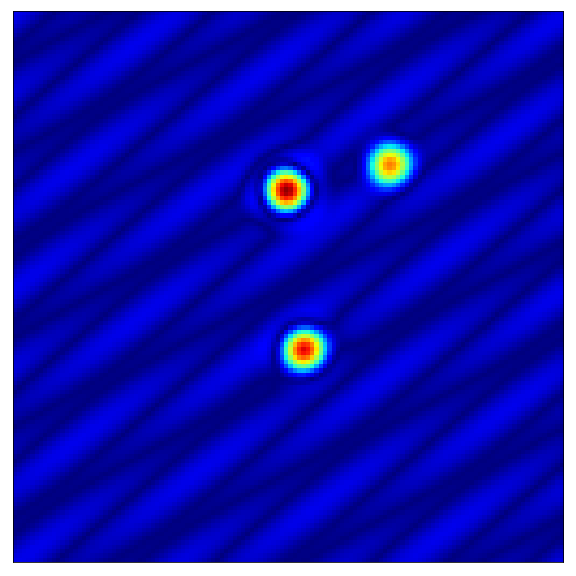}
    \hspace{0.1in}
    \includegraphics[height=2.1in]{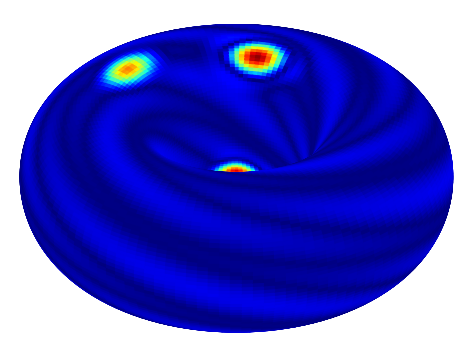}
     
    \caption{A snapshot of the random system in 2D (left) and on a 3D torus (right).}
    \label{fig:torus}
\end{figure}

\section{Numerical Experiments} \label{sec:num_exp}
In practice, it is often the case that for spatiotemporal data and systems, the underlying low-rank dynamic model, the disturbance, and the noise are not known. 
In this case, we would fit an estimated model representation from data via DMD, and approximate disturbance and noise covariances either from data or through hyper-parameter tuning. 
Here, we look at two examples: (i) the Kuramoto Sivashinsky (KS) system, and (ii) the Sea Surface Temperature (SST) dataset from NOAA~\cite{noaasst,reynolds2002improved}. 
We study the performance when we use a DMD approximated model for Kalman filter estimation and sensor path planning by our greedy algorithm. 
In both examples, we fit a linear DMD model with a chosen low rank.
The dynamics matrix $\mathbf{\Lambda}$ is diagonal with DMD eigenvalues and the basis $\mathbf{\Psi}$ consists of the DMD modes. 
We further add a white measurement noise with variance $\mathbf{R} = \mathbf{I}$ to the data, and we set the system disturbance to be $\mathbf{Q} = q\mathbf{I}$ where the uniform variance $q$ is a hyperparameter tuned by experiment.

These two examples are representative in different aspects. 
The KS system is known for its chaotic behavior. Therefore, a linear representation of the system is extremely difficult. Additionally, the modes from the linear approximated model are mostly local since linear correlation between locations is small, so that full observability is hard to achieve with few fixed sensors. We show in Section \ref{sec:ks} that mobile sensors can perform particularly well comparing to fixed sensors by reaching more locations and capturing more local structures. 

The SST dataset from NOAA contains weekly mean optimum interpolated sea surface temperature measurements from global satellite data. The dataset can be well approximated by a linear model and most modes in the approximated system are global so that observability is easily achieved with even just one stationary sensor. We show then in Section \ref{sec:sst} mobile sensors further accelerate the convergence of error.

\subsection{Kuramoto Sivashinsky System} \label{sec:ks}

The KS system is given by the equation $u_t + uu_x + u_{xx} + u_{xxxx} = 0$.
We consider the numerical solution of the system on a spatial grid of size 2048 over $x\in[0, 2\pi]$. The initial condition is randomly generated over a standard normal distribution. 
With a random initial condition, we numerically solve the KS equation and collect data on the time interval $t \in [0,10]$ with a time step of $dt=10^{-4}$.
We first perform singular value decomposition (SVD) to find a low rank representation of the data. 
The first 100 singular values capture 99.99\% of the energy, so we estimate a low-dimensional linear representation of the system by fitting a standard dynamic mode decomposition (DMD) model~\cite{Tu2014jcd} with an SVD rank of 100.

\begin{figure}[t]
    \centering
    \begin{overpic}[abs,height=2in]{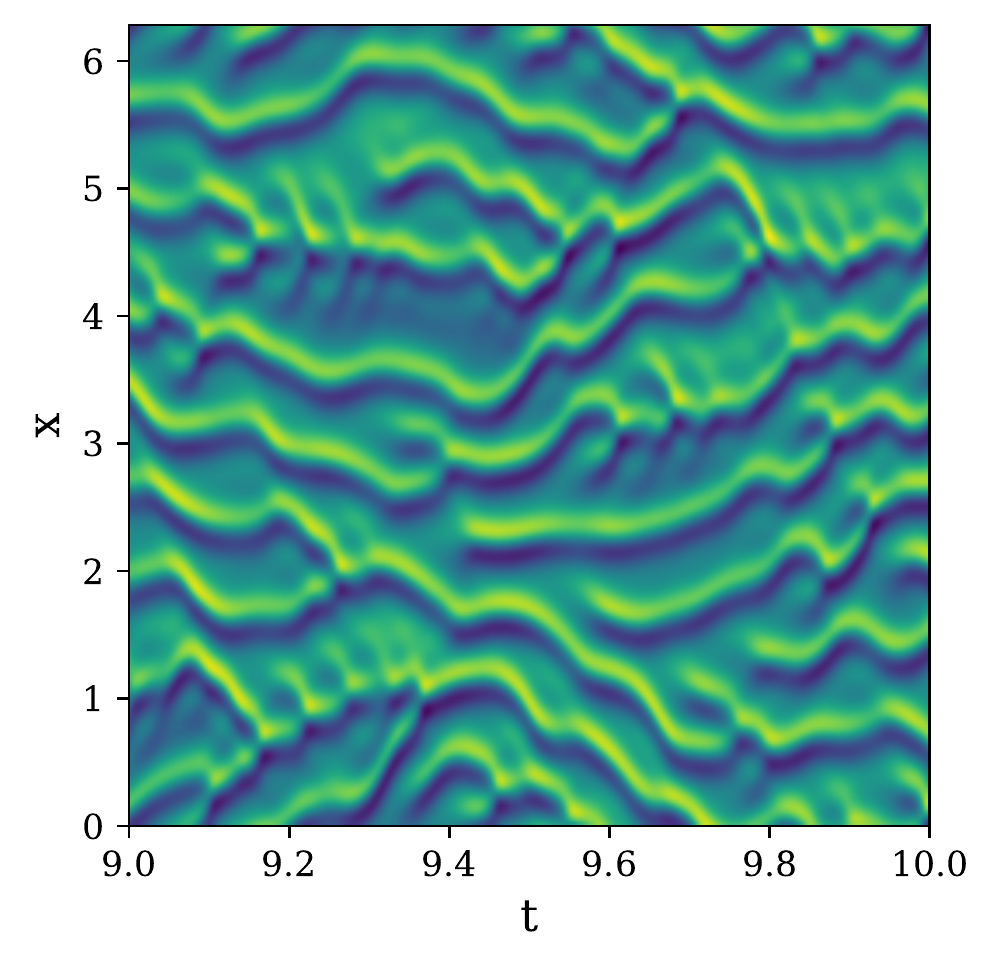}
        \put(-15,130){(a)}
    \end{overpic}
    \hspace{.3in}
    \begin{overpic}[abs,height=2in]{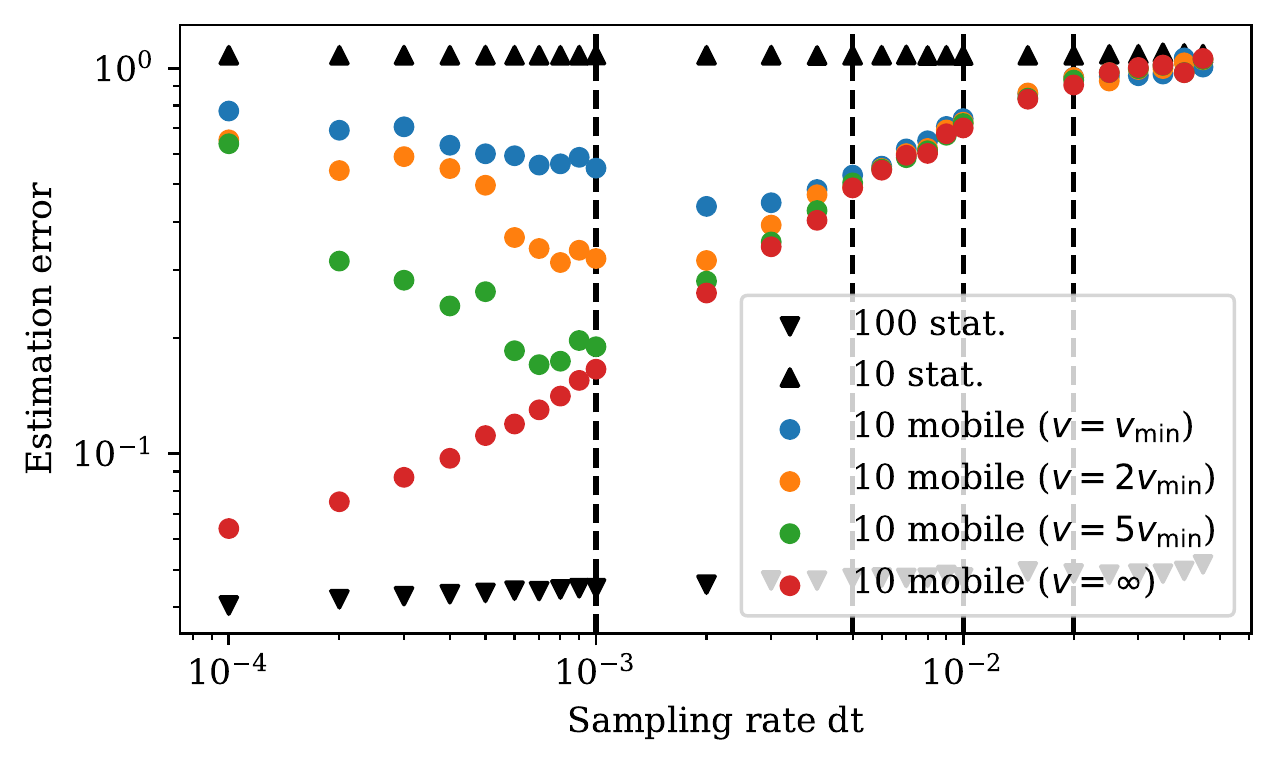}
        \put(-15,130){(b)}
    \end{overpic}
    \vfill
    \hspace{.2in}
    \begin{overpic}[width=.23\textwidth]{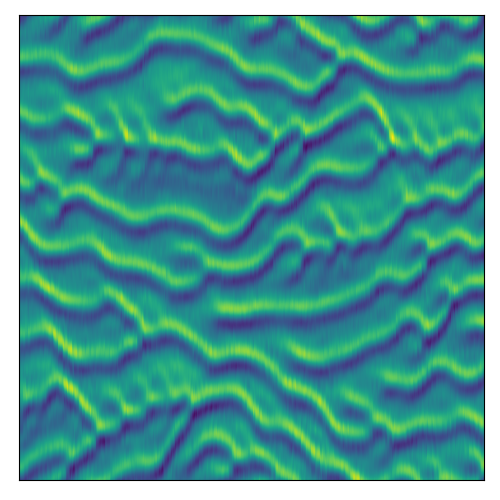}
        \put(5,5){\color{white}\boldmath$dt=0.001$}
        \put(-15,90){(c)}
    \end{overpic}
    \begin{overpic}[width=.23\textwidth]{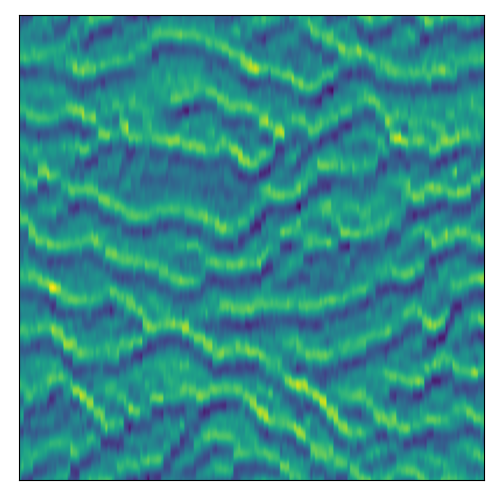}
        \put(5,5){\color{white}\boldmath$dt=0.005$}
    \end{overpic}
    \begin{overpic}[width=.23\textwidth]{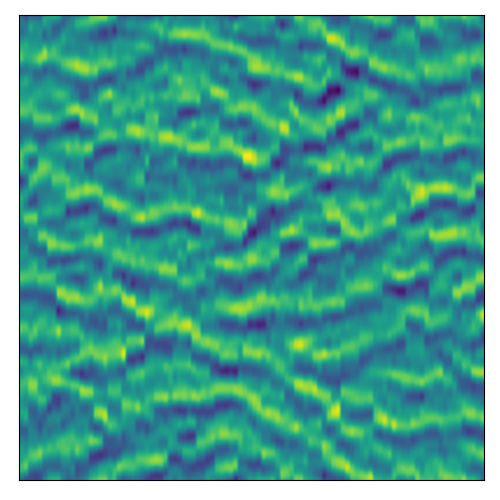}
        \put(5,5){\color{white}\boldmath$dt=0.01$}
    \end{overpic}
    \begin{overpic}[width=.23\textwidth]{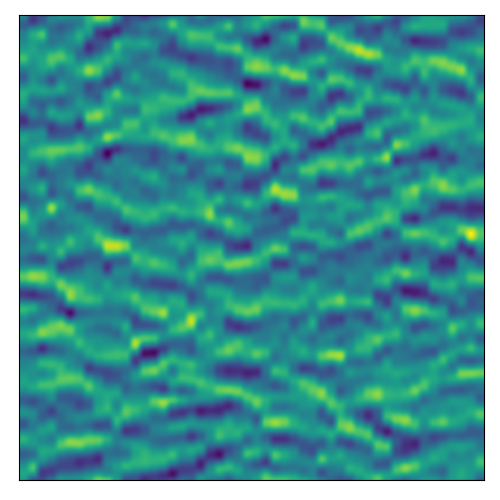}
        \put(5,5){\color{white}\boldmath$dt=0.02$}
    \end{overpic}
     
    \caption{(a) True spatiotemporal dynamics of the KS system in $T \in [9,10]$; (b) Bode plot of estimation error against sampling rate; (c) Estimated x-t plot by 10 mobile sensors with sampling rate $dt = 0.001, 0.005, 0.01, 0.02$ (corresponding with the dashed vertial lines on the bode plot).}
    \label{fig:ks}
\end{figure}

Because of the chaotic nature of the system, accurate estimation is not possible with a limited number of 10 sparse fixed sensors.  Indeed, we need 100 fixed sensors, equivalent to the full rank of our approximated linear system, to effectively estimate the system (Figure \ref{fig:ks}).
Additionally, we see that there is no significant improvement in performance with increasing sampling rate using fixed sensors since more frequent measurements at the same locations add little information of the unobserved states. 

On the other hand, mobile sensors can move to measure different locations and gain more information of the entire state space. 
With fast enough sensor speed, 10 moving sensors can achieve a significantly improved estimation comparing to 10 fixed sensors by increasing the sampling rate (Figure \ref{fig:ks}). 
The improvement in performance is limited by the sensor speed. 
We set the minimum  speed in this example to be $v_{\min}=\frac{2\pi}{2048}*10^4$ so that the sensor is able to move to its left and right neighbor at a discrete sampling time step of $10^{-4}$. 
With higher sensor speed, sensors can make observations over a wider spatial range, thus giving better estimation. 
As $v \to \infty$, the performance of 10 moving sensors approaches that of 100 fixed sensors with fast enough sampling rate.

\begin{figure}[t]
    \centering
    \includegraphics[width=.6\textwidth]{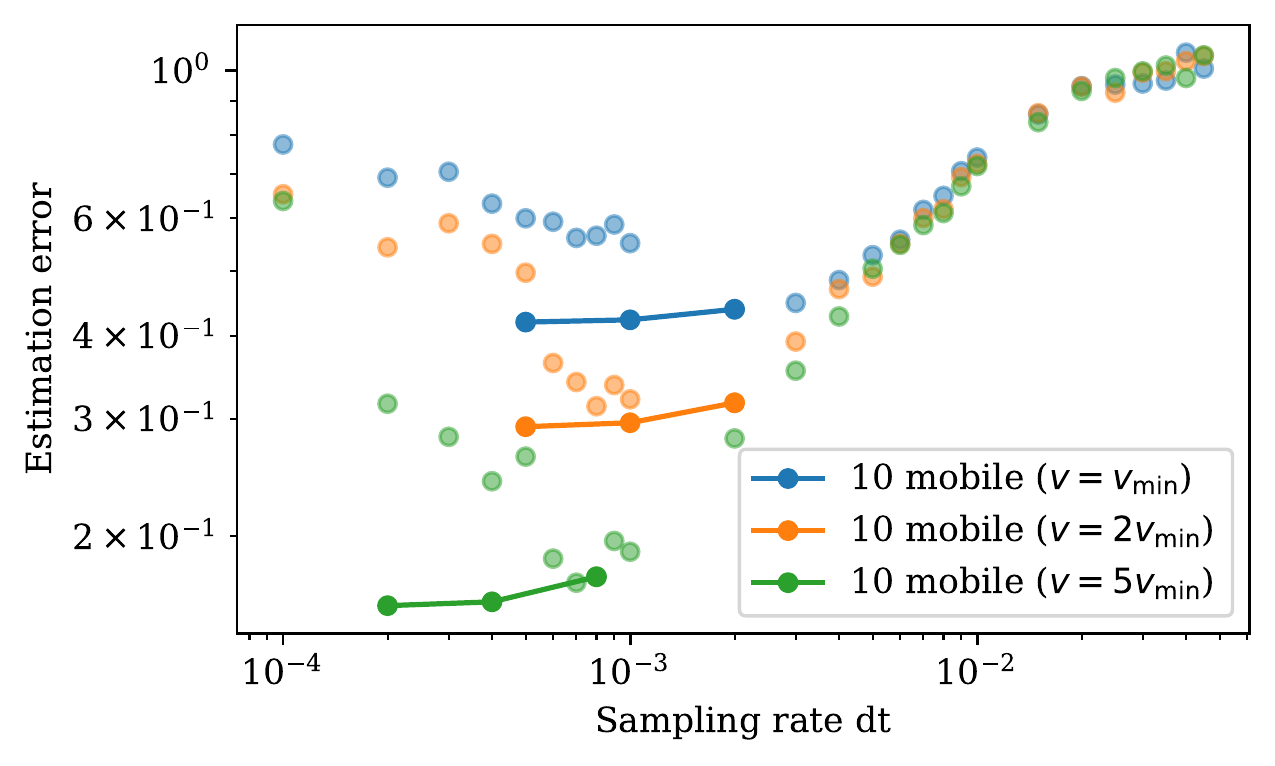}
    \caption{Estimation error against sampling rate with multiscale expansion plotted in full transparency connected by line. }
    \label{fig:ks_multi}
\end{figure}

Due to the greedy nature of our algorithm, it selects based on the immediate reward at the next time step and cannot look ahead. 
When the sampling is sufficiently fast, the greedy algorithm makes a decision based on the closest neighbors of the current location. 
Such a decision is not informative, and the trajectory planned fails to have a good performance. 
One way to reduce the greediness of the algorithm is to perform a multiscale path completion. 
We start by finding a trajectory at a slower sampling rate.
Then, we gradually decrease the time step and apply the same path planning algorithm, except using the previously found trajectory as guidance and filling in the gap to construct a more complete trajectory at the faster sampling rate. 
We apply this multiscale expansion procedure on the KS example, initiated at the sampling rate with the smallest error, and expand to faster sampling rate based on that path.
We see the performance is no longer worse with a fast sampling rate in the KS experiments, but it flattens and reaches a limit determined by the sensor speed (Figure \ref{fig:ks_multi}).

\subsection{Sea Surface Temperature} \label{sec:sst}
The Sea Surface Temperature (SST) dataset contains weekly collection from satellite data of sea surface temperature measurements on the 1 degree latitude by 1 degree longitude (180 by 360) global grid from 1990 to the present. 
We fit a standard DMD model with a low dimension of 10.

\begin{figure}[t]
    \centering
    \hspace{.04\textwidth}
    \begin{overpic}[width=.45\textwidth]{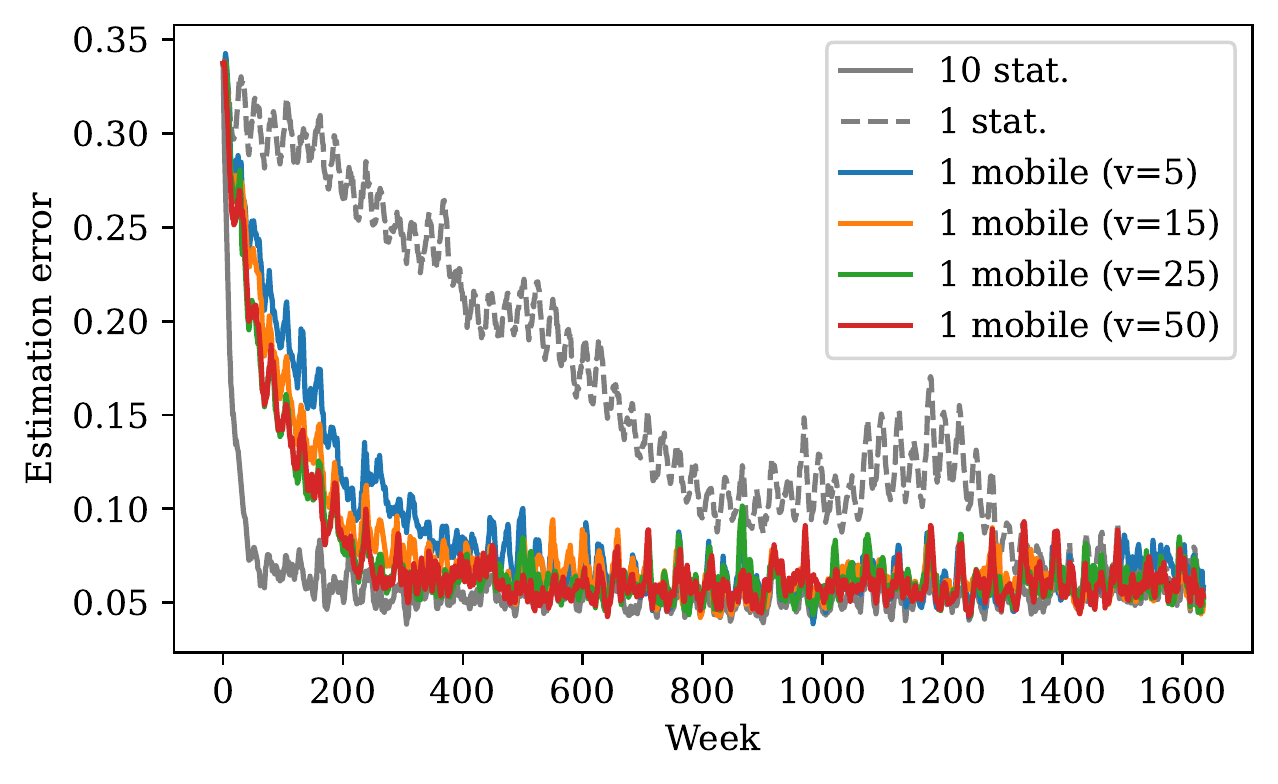}
        \put(-7,56){(a)}
    \end{overpic}
    \hspace{.04\textwidth}
    \begin{overpic}[width=.45\textwidth]{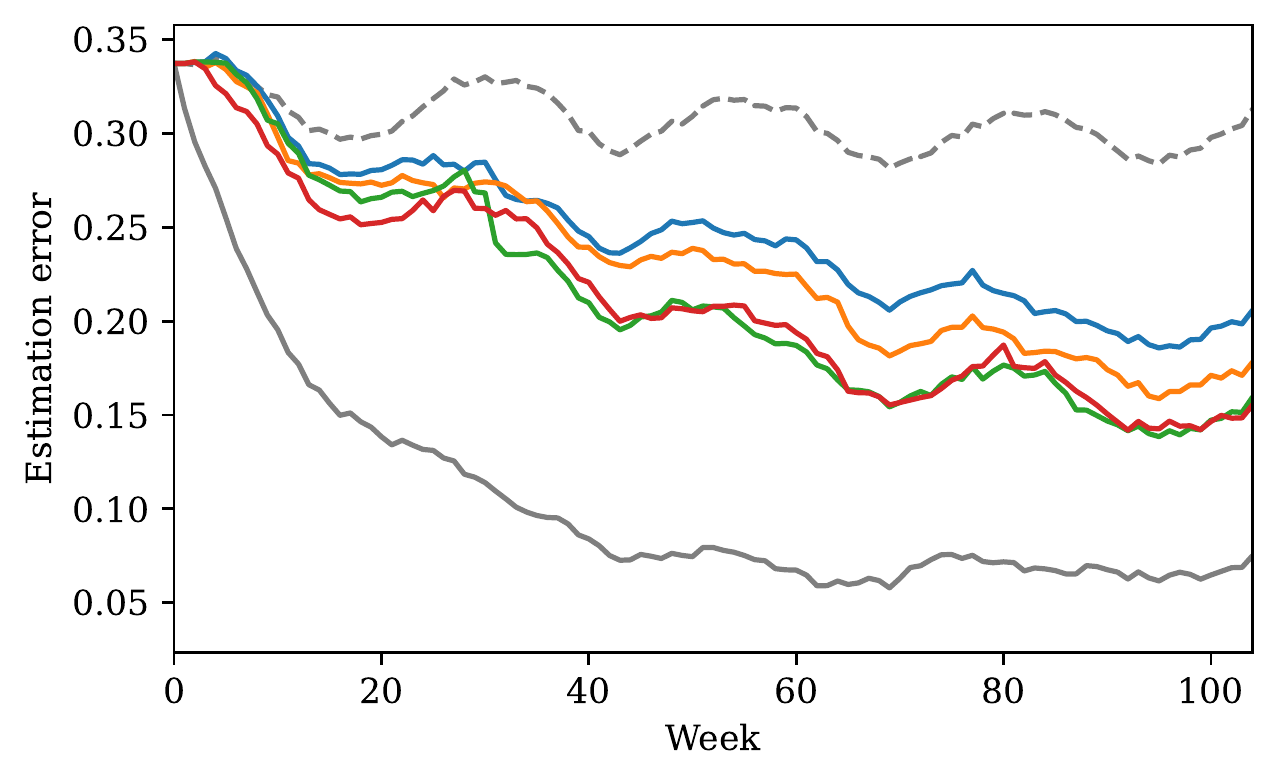}
        \put(-7,56){(b)}
    \end{overpic}
    
    \caption{Estimation error over (a) all time; (b) first two year (104 weekly measurements). }
    \label{fig:sst}
\end{figure}

First, we see that Kalman filter estimation using one stationary sensor is comparable to that of ten sensors (Figure \ref{fig:sst}). 
This verifies that the approximated linear model contains mostly global features that can be observed well with very few sensors. 
However, with a bad initial estimation, KF with 10 fixed sensors converges to low error very quickly (below 0.1 within one year), while it takes 1 sensor almost 28 years to reach a comparable error. 

Next, we consider estimating the SST data with one mobile sensor. 
Since the limiting estimation is already good with one stationary sensor, we show the improvements in the KF convergence speed for mobile sensors.
The DMD eigenvalues show the max frequency of the system to be around a half-year, so we choose the period of mobile sensor trajectory to be about a quarter-year (14 weeks). 
Since the globe has continental land as obstructions, we need to ensure the planned trajectory does not cross any land as the sensor moves in water. 
We build a connectivity graph and adjacency matrix for the candidate selection step in our algorithm instead of a simple Euclidean distance function. 

\begin{figure}[t]
    \centering
    \begin{overpic}[width=.9\textwidth]{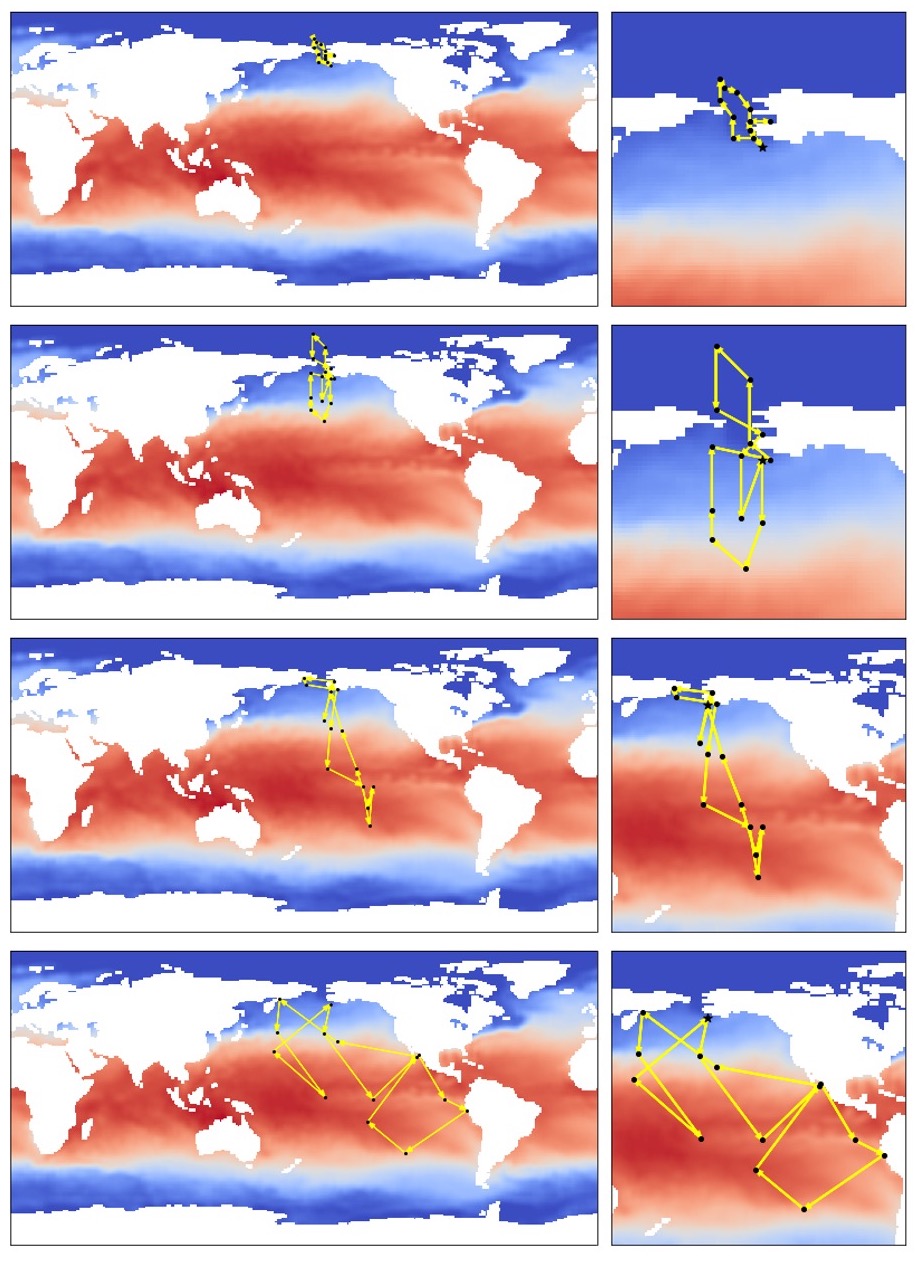} 
        \put(49,3){\color{black} \boldmath$v=50$}
        \put(49,27.5){\color{black} \boldmath$v=25$}    
        \put(49,52.5){\color{black} \boldmath$v=15$}    
        \put(49,77){\color{black} \boldmath$v=5$}    
    \end{overpic}
    \caption{Planned sensor trajectory with a cycle period of 14 weeks, where the movement speed is limited to 5, 15, 25, 50 spatial units (1 degree of latitude or longitude). Zoomed in map on the right.}
    \label{fig:sst_path}
\end{figure}

\begin{figure}[t]
    \centering
    \includegraphics[width=.7\textwidth]{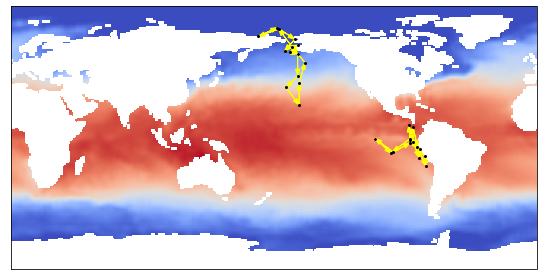}
    \caption{Planned sensor trajectory for 2 sensors with a cycle period of 14 weeks, with sensor speed limit at 15 spatial units (1 degree of latitude or longitude).}
    \label{fig:sst_path2}
\end{figure}

Figure \ref{fig:sst} shows the results of one mobile sensor with different sensor speed limits. 
Mobile sensor estimation indeed produces much faster convergence compared to a stationary sensor. 
As the speed limit increases, the sensor can move to farther locations with better observability, further improving the convergence of estimation. 
Figure \ref{fig:sst_path} shows the paths of the sensor.
When $v$ is small, the initial location plays an important role since the trajectory does not move far from it. 
The first location picked by the algorithm is close to Alaska, so the first two trajectories with low speed center around the North Pacific and the Arctic Ocean. 
As $v$ increases, the sensor explores the equator and south hemisphere regions, especially the El Nino regions around the equatorial Pacific, which is an important local feature.
Figure \ref{fig:sst_path2} shows the planned trajectory using 2 mobile sensors.

\begin{figure}[t]
    \centering
    \begin{overpic}[width=.7\textwidth]{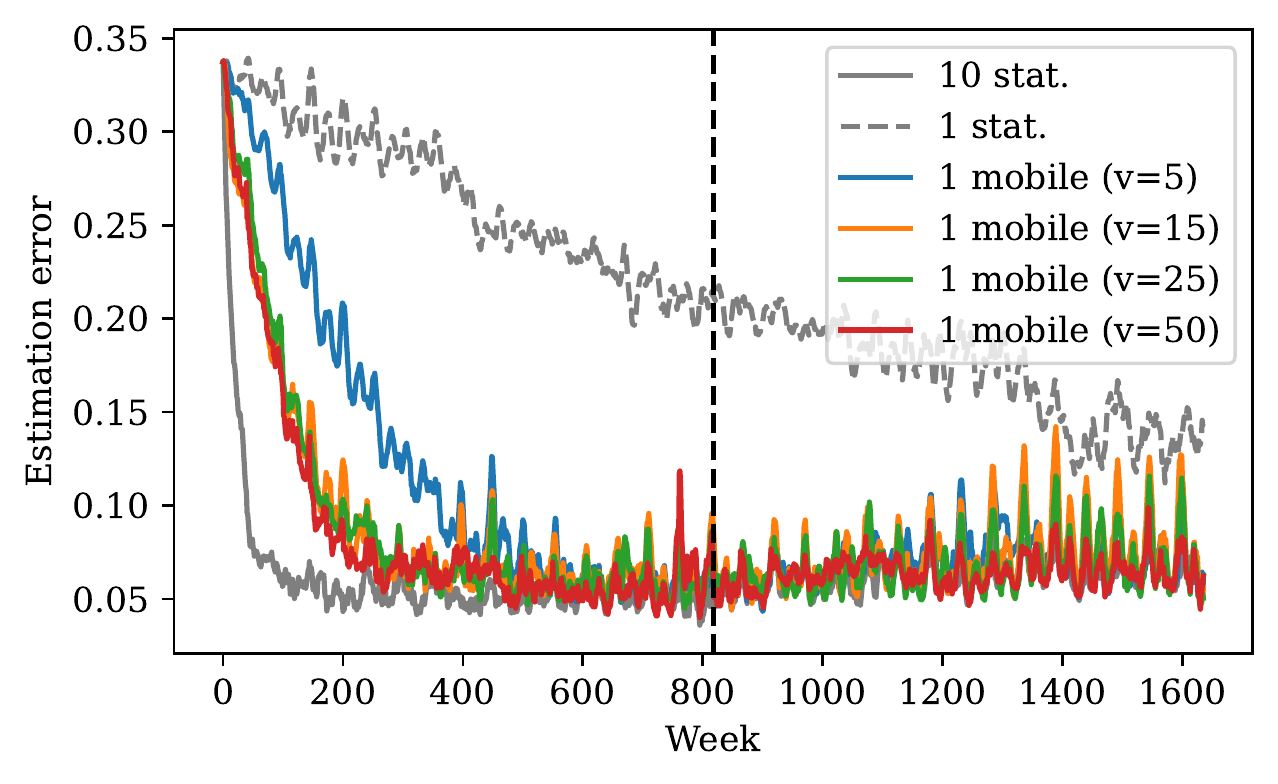}
        \put(38,60){train $\leftarrow$}
        \put(58,60){$\rightarrow$ test}
    \end{overpic}
    
    \caption{Estimation error using approximated DMD model trained on the first half of the data. }
    \label{fig:sst_train}
\end{figure}

KF estimation convergence also matters when the underlying dynamics is nonstationary and changes over time. 
If the estimation does not reach a meaningful error in time, the shifting dynamics will further slow down the convergence and increase the limiting error.
To show this, we instead fit a DMD model using only the first half of the SST data and use it as the approximated linear model for KF estimation. 
In this case, the fitted linear model is representative and relevant only in the first training half, and does not reflect any possible changes in the data dynamics afterwards.
Then, one stationary sensor performs significantly worse due to slow convergence (Figure \ref{fig:sst_train}). 
On the other hand, the error from one mobile sensor converges fast enough within the training period so that in the second half the error is still relatively low. 
Therefore, fast convergence with mobile sensors ensures a fast adjustment in estimation when the dynamics change in time.

\clearpage

\section{Conclusion and Future Work}
In this work, we developed a mathematical strategy for planning a periodic trajectory for limited mobile sensors to estimate a spatiotemporal system using Kalman filter estimation. 
We examine the system observability as a metric that influences the estimation performance in terms of the limiting squared error as well as the convergence rate.
We consider an objective to minimize the condition number of the discrete observability matrix along the trajectory and formulate it as a submatrix selection problem. 
We then propose a time-forwarding greedy algorithm that selects sensor locations along the trajectory using the same rules as QRcp and GappyPOD+E from a carefully chosen candidate subset.

The experiments show that the method is able to plan a trajectory that locates the local features and improves the estimation performance. 
In these experiments, we explore Kalman filter design factors  and their impact on estimation as they relate to the three important timescales: the Nyquist rate of the underlying dynamics, the rate of sampling, and the velocity of the sensors.   We find that mobile sensors are especially beneficial for a complex, non-linear system to capture local features in an approximated linear model, without deploying a large amount of sensors. 
We also see an improvement in estimation convergence rate using mobile sensors, which more rapidly reaches an accurate estimation. 

In future work, a weighted cost function can be added to the objective to better incorporate different costs to the path planning. 
Sensor speed can be formulated as a cost instead of a hard constraint imposed in the selection process. 
Furthermore, energy consumption caused by sensor movement can also be included to plan a trajectory that is also more energy efficient. 
For example, as in the flow field applications, we can consider the flow field information and the energy cost associated with it as the sensor moves with or against the flow. 
Incorporating the background flow velocity in the set of possible next locations is an important future extension of this work.  
We can also refer to many different sensor control laws as cost constraints for incorporating other tasks such as simultaneous structure tracking.

For multi-sensor planning, it will be interesting to consider different asynchronous periodic trajectory for each individual sensor instead of all having the same period. This will be particularly useful for multiscale systems so that each sensor can be responsible for estimating features of different timescales.
As shown in the numerical experiments, the performance of Kalman filter estimation fundamentally depends on the accurate modeling of the linear system and the correct choices of hyper-parameters. 
Better data-driven linear system identification can be explored. 
An alternating model fitting and estimation approach can be explored to update both the model and the sensor trajectory continuously to achieve even better performance.

\section*{Acknowledgments}
The authors acknowledge support from the National Science Foundation AI Institute in Dynamic Systems (grant number 2112085). 
SLB acknowledges support from the Air Force Office of Scientific Research (FA9550-21-1-0178). 
JNK acknowledges support from the Air Force Office of Scientific Research (FA9550-19-1-0011).

\printbibliography

@article{krause2008near,
  title={Near-optimal sensor placements in Gaussian processes: Theory, efficient algorithms and empirical studies.},
  author={Krause, Andreas and Singh, Ajit and Guestrin, Carlos},
  journal={Journal of Machine Learning Research},
  volume={9},
  number={2},
  year={2008}
}

@article{caselton1984optimal,
  title={Optimal monitoring network designs},
  author={Caselton, William F and Zidek, James V},
  journal={Statistics \& Probability Letters},
  volume={2},
  number={4},
  pages={223--227},
  year={1984},
  publisher={Elsevier}
}

@article{li2021efficient,
  title={Efficient Sensor Placement for Signal Reconstruction Based on Recursive Methods},
  author={Li, Bangjun and Liu, Haoran and Wang, Ruzhu},
  journal={IEEE Transactions on Signal Processing},
  volume={69},
  pages={1885--1898},
  year={2021},
  publisher={IEEE}
}

@article{Tropp:2004,
	author = {J. A. Tropp},
	journal = {IEEE Transactions on Information Theory},
	number = {10},
	pages = {2231--2242},
	title = {Greed is good: Algorithmic results for sparse approximation},
	volume = {50},
	year = {2004}}

@article{Tropp:2006b,
	author = {Tropp, J. A. and Gilbert, A. C. and Strauss, M. J.},
	journal = {Signal Processing},
	number = {3},
	pages = {572--588},
	title = {Algorithms for simultaneous sparse approximation. Part I: Greedy pursuit},
	volume = {86},
	year = {2006}}

@article{manohar2018data,
  title={Data-driven sparse sensor placement for reconstruction: Demonstrating the benefits of exploiting known patterns},
  author={Manohar, Krithika and Brunton, Bingni W and Kutz, J Nathan and Brunton, Steven L},
  journal={IEEE Control Systems Magazine},
  volume={38},
  number={3},
  pages={63--86},
  year={2018},
  publisher={IEEE}
}

@article{clark2018greedy,
  title={Greedy sensor placement with cost constraints},
  author={Clark, Emily and Askham, Travis and Brunton, Steven L and Kutz, J Nathan},
  journal={IEEE Sensors Journal},
  volume={19},
  number={7},
  pages={2642--2656},
  year={2018},
  publisher={IEEE}
}

@article{peherstorfer_stability_2020,
	title = {Stability of Discrete Empirical Interpolation and Gappy Proper Orthogonal Decomposition with Randomized and Deterministic Sampling Points},
	volume = {42},
	issn = {1064-8275, 1095-7197},
	doi = {10.1137/19M1307391},
	pages = {A2837--A2864},
	number = {5},
	journaltitle = {{SIAM} Journal on Scientific Computing},
	shortjournal = {{SIAM} J. Sci. Comput.},
	author = {Peherstorfer, Benjamin and Drmač, Zlatko and Gugercin, Serkan},
	date = {2020-01},
	langid = {english}
}

@article{astrid2008missing,
  title={Missing point estimation in models described by proper orthogonal decomposition},
  author={Astrid, Patricia and Weiland, Siep and Willcox, Karen and Backx, Ton},
  journal={IEEE Transactions on Automatic Control},
  volume={53},
  number={10},
  pages={2237--2251},
  year={2008},
  publisher={IEEE}
}

@article{saito2021determinant,
  title={Determinant-based fast greedy sensor selection algorithm},
  author={Saito, Yuji and Nonomura, Taku and Yamada, Keigo and Nakai, Kumi and Nagata, Takayuki and Asai, Keisuke and Sasaki, Yasuo and Tsubakino, Daisuke},
  journal={IEEE Access},
  volume={9},
  pages={68535--68551},
  year={2021},
  publisher={IEEE}
}

@article{drmac2016new,
  title={A new selection operator for the discrete empirical interpolation method---improved a priori error bound and extensions},
  author={Drmac, Zlatko and Gugercin, Serkan},
  journal={SIAM Journal on Scientific Computing},
  volume={38},
  number={2},
  pages={A631--A648},
  year={2016},
  publisher={SIAM}
}

@article{manohar_optimal_2021,
	title = {Optimal Sensor and Actuator Selection using Balanced Model Reduction},
	issn = {1558-2523},
	doi = {10.1109/TAC.2021.3082502},
	pages = {1--1},
	journaltitle = {{IEEE} Transactions on Automatic Control},
	author = {Manohar, Krithika and Kutz, J. Nathan and Brunton, Steven L.},
	date = {2021}
}

@article{williams2022data,
  title={Data-driven sensor placement with shallow decoder networks},
  author={Williams, Jan and Zahn, Olivia and Kutz, J Nathan},
  journal={arXiv preprint arXiv:2202.05330},
  year={2022}
}

@ARTICLE{clark2020cost,
  author={Clark, Emily and Kutz, J. Nathan and Brunton, Steven L.},
  journal={IEEE Sensors Journal}, 
  title={Sensor Selection With Cost Constraints for Dynamically Relevant Bases}, 
  year={2020},
  volume={20},
  number={19},
  pages={11674-11687},
  doi={10.1109/JSEN.2020.2997298}}

@inproceedings{tzoumas_sensor_2016,
	title = {Sensor placement for optimal Kalman filtering: Fundamental limits, submodularity, and algorithms},
	doi = {10.1109/ACC.2016.7524914},
	shorttitle = {Sensor placement for optimal Kalman filtering},
	eventtitle = {2016 American Control Conference ({ACC})},
	pages = {191--196},
	booktitle = {2016 American Control Conference ({ACC})},
	author = {Tzoumas, V. and Jadbabaie, A. and Pappas, G. J.},
	date = {2016-07}
}

@inproceedings{ye_complexity_2018,
	title = {On the Complexity and Approximability of Optimal Sensor Selection for Kalman Filtering},
	doi = {10.23919/ACC.2018.8431016},
	eventtitle = {2018 Annual American Control Conference ({ACC})},
	pages = {5049--5054},
	booktitle = {2018 Annual American Control Conference ({ACC})},
	author = {Ye, Lintao and Roy, Sandip and Sundaram, Shreyas},
	date = {2018-06},
}

@inproceedings{dhingra_admm_2014,
	title = {An {ADMM} algorithm for optimal sensor and actuator selection},
	doi = {10.1109/CDC.2014.7040017},
	eventtitle = {53rd {IEEE} Conference on Decision and Control},
	pages = {4039--4044},
	booktitle = {53rd {IEEE} Conference on Decision and Control},
	author = {Dhingra, Neil K. and Jovanović, Mihailo R. and Luo, Zhi-Quan},
	date = {2014-12}
}

@inproceedings{lan_planning_2013,
	title = {Planning periodic persistent monitoring trajectories for sensing robots in Gaussian Random Fields},
	doi = {10.1109/ICRA.2013.6630905},
	eventtitle = {2013 {IEEE} International Conference on Robotics and Automation},
	pages = {2415--2420},
	booktitle = {2013 {IEEE} International Conference on Robotics and Automation},
	author = {Lan, Xiaodong and Schwager, Mac},
	date = {2013-05}
}

@article{everson1995karhunen,
  title={Karhunen--Loeve procedure for gappy data},
  author={Everson, Richard and Sirovich, Lawrence},
  journal={JOSA A},
  volume={12},
  number={8},
  pages={1657--1664},
  year={1995},
  publisher={Optica Publishing Group}
}

@article{erichson2020shallow,
  title={Shallow neural networks for fluid flow reconstruction with limited sensors},
  author={Erichson, N Benjamin and Mathelin, Lionel and Yao, Zhewei and Brunton, Steven L and Mahoney, Michael W and Kutz, J Nathan},
  journal={Proceedings of the Royal Society A},
  volume={476},
  number={2238},
  pages={20200097},
  year={2020},
  publisher={The Royal Society Publishing}
}

@book{trefethen1997numerical,
  title={Numerical linear algebra},
  author={Trefethen, Lloyd N and Bau III, David},
  volume={50},
  year={1997},
  publisher={Siam}
}

@INPROCEEDINGS{zhang_optimal_nodate,  
    author={Zhang, Wei and Vitus, Michael P. and Hu, Jianghai and Abate, Alessandro and Tomlin, Claire J.},  
    booktitle={49th IEEE Conference on Decision and Control (CDC)},   
    title={On the optimal solutions of the infinite-horizon linear sensor scheduling problem},   
    year={2010},  
    volume={},  
    number={},  
    pages={396-401},  
    doi={10.1109/CDC.2010.5717163}
}

@article{zhao_optimal_2014,
	title = {On the Optimal Solutions of the Infinite-Horizon Linear Sensor Scheduling Problem},
	volume = {59},
	issn = {1558-2523},
	doi = {10.1109/TAC.2014.2314222},
	pages = {2825--2830},
	number = {10},
	journaltitle = {{IEEE} Transactions on Automatic Control},
	author = {Zhao, Lin and Zhang, Wei and Hu, Jianghai and Abate, Alessandro and Tomlin, Claire J.},
	date = {2014-10}
}

@article{liu_optimal_2014,
	title = {Optimal Periodic Sensor Scheduling in Networks of Dynamical Systems},
	volume = {62},
	issn = {1941-0476},
	doi = {10.1109/TSP.2014.2320455},
	pages = {3055--3068},
	number = {12},
	journaltitle = {{IEEE} Transactions on Signal Processing},
	author = {Liu, Sijia and Fardad, Makan and Masazade, Engin and Varshney, Pramod K.},
	date = {2014-06}
}

@article{lan_rapidly_2016,
	title = {Rapidly Exploring Random Cycles: Persistent Estimation of Spatiotemporal Fields With Multiple Sensing Robots},
	volume = {32},
	issn = {1941-0468},
	doi = {10.1109/TRO.2016.2596772},
	shorttitle = {Rapidly Exploring Random Cycles},
	pages = {1230--1244},
	number = {5},
	journaltitle = {{IEEE} Transactions on Robotics},
	author = {Lan, Xiaodong and Schwager, Mac},
	date = {2016-10},

}

@article{chen_deep_2020,
	title = {Deep Reinforced Learning Tree for Spatiotemporal Monitoring With Mobile Robotic Wireless Sensor Networks},
	volume = {50},
	issn = {2168-2232},
	doi = {10.1109/TSMC.2019.2920390},
	pages = {4197--4211},
	number = {11},
	journaltitle = {{IEEE} Transactions on Systems, Man, and Cybernetics: Systems},
	author = {Chen, Jiahong and Shu, Tongxin and Li, Teng and de Silva, Clarence W.},
	date = {2020-11},
	
}

@article{mo_infinite-horizon_2014,
	title = {On infinite-horizon sensor scheduling},
	volume = {67},
	issn = {01676911},
	doi = {10.1016/j.sysconle.2014.02.002},
	pages = {65--70},
	journaltitle = {Systems \& Control Letters},
	shortjournal = {Systems \& Control Letters},
	author = {Mo, Y. and Garone, E. and Sinopoli, B.},
	date = {2014-05},
	langid = {english},

}

@article{shi_approximate_2013,
	title = {Approximate optimal periodic scheduling of multiple sensors with constraints},
	volume = {49},
	issn = {00051098},
	doi = {10.1016/j.automatica.2013.01.024},
	pages = {993--1000},
	number = {4},
	journaltitle = {Automatica},
	shortjournal = {Automatica},
	author = {Shi, Dawei and Chen, Tongwen},
	date = {2013-04},
	langid = {english},

}

@article{zhang_sensor_2017,
	title = {Sensor selection for Kalman filtering of linear dynamical systems: Complexity, limitations and greedy algorithms},
	volume = {78},
	issn = {00051098},
	doi = {10.1016/j.automatica.2016.12.025},
	shorttitle = {Sensor selection for Kalman filtering of linear dynamical systems},
	pages = {202--210},
	journaltitle = {Automatica},
	shortjournal = {Automatica},
	author = {Zhang, Haotian and Ayoub, Raid and Sundaram, Shreyas},
	date = {2017-04},
	langid = {english}
}

@article{asghar_complete_2017,
	title = {A complete greedy algorithm for infinite-horizon sensor scheduling},
	volume = {81},
	issn = {00051098},
	doi = {10.1016/j.automatica.2017.04.018},
	pages = {335--341},
	journaltitle = {Automatica},
	shortjournal = {Automatica},
	author = {Asghar, Ahmad Bilal and Jawaid, Syed Talha and Smith, Stephen L.},
	date = {2017-07},
	langid = {english},
	
}

@thesis{ilkturk_observability_nodate,
	location = {United States -- Arizona},
	title = {Observability Methods in Sensor Scheduling},
	rights = {Database copyright {ProQuest} {LLC}; {ProQuest} does not claim copyright in the individual underlying works.},
	pagetotal = {86},
	institution = {Arizona State University},
	type = {phdthesis},
	author = {Ilkturk, Utku},
	year = {2015},
	
}

@article{rafieisakhaei_use_2017,
	title = {On the Use of the Observability Gramian for Partially Observed Robotic Path Planning Problems},
	doi = {10.1109/CDC.2017.8263868},
	pages = {1523--1528},
	journaltitle = {2017 {IEEE} 56th Annual Conference on Decision and Control ({CDC})},
	author = {Rafieisakhaei, Mohammadhussein and Chakravorty, Suman and Kumar, P. R.},
	date = {2017-12},
	eprinttype = {arxiv},
	eprint = {1801.09877},
	
}

@article{chamon_approximate_2021-1,
	title = {Approximate Supermodularity of Kalman Filter Sensor Selection},
	volume = {66},
	issn = {1558-2523},
	doi = {10.1109/TAC.2020.2973774},
	pages = {49--63},
	number = {1},
	journaltitle = {{IEEE} Transactions on Automatic Control},
	author = {Chamon, Luiz F. O. and Pappas, George J. and Ribeiro, Alejandro},
	date = {2021-01},
	
}

@article{komaroff1992upper,
  title={Upper bounds for the solution of the discrete Riccati equation},
  author={Komaroff, N},
  journal={IEEE transactions on automatic control},
  volume={37},
  number={9},
  pages={1370--1373},
  year={1992},
  publisher={IEEE}
}

@article{komaroff1992lower,
  title={Lower summation bounds for the discrete Riccati and Lyapunov equations},
  author={Komaroff, N and Shahian, B},
  journal={IEEE Transactions on Automatic Control},
  volume={37},
  number={7},
  pages={1078--1080},
  year={1992},
  publisher={IEEE}
}

@incollection{bittanti1991periodic,
  title={The periodic Riccati equation},
  author={Bittanti, Sergio and Colaneri, Patrizio and Nicolao, Giuseppe De},
  booktitle={The Riccati Equation},
  pages={127--162},
  year={1991},
  publisher={Springer}
}

@article{dai2011eigenvalue,
  title={On eigenvalue bounds and iteration methods for discrete algebraic Riccati equations},
  author={Dai, Hua and Bai, Zhong-Zhi},
  journal={Journal of Computational Mathematics},
  pages={341--366},
  year={2011},
  publisher={JSTOR}
}

@article{kwon1996bounds,
  title={Bounds in algebraic Riccati and Lyapunov equations: a survey and some new results},
  author={Kwon, Wook Hyun and Moon, Young Soo and Ahn, Sang Chul},
  journal={International Journal of Control},
  volume={64},
  number={3},
  pages={377--389},
  year={1996},
  publisher={Taylor \& Francis}
}

@article{leonard2007collective,
  title={Collective motion, sensor networks, and ocean sampling},
  author={Leonard, Naomi Ehrich and Paley, Derek A and Lekien, Francois and Sepulchre, Rodolphe and Fratantoni, David M and Davis, Russ E},
  journal={Proceedings of the IEEE},
  volume={95},
  number={1},
  pages={48--74},
  year={2007},
  publisher={IEEE}
}

@article{paley2020mobile,
  title={Mobile sensor networks and control: Adaptive sampling of spatiotemporal processes},
  author={Paley, Derek A and Wolek, Artur},
  journal={Annual Review of Control, Robotics, and Autonomous Systems},
  volume={3},
  pages={91--114},
  year={2020},
  publisher={Annual Reviews}
}

@article{devries2013observability,
  title={Observability-based optimization of coordinated sampling trajectories for recursive estimation of a strong, spatially varying flowfield},
  author={DeVries, Levi and Majumdar, Sharanya J and Paley, Derek A},
  journal={Journal of intelligent \& robotic systems},
  volume={70},
  number={1},
  pages={527--544},
  year={2013},
  publisher={Springer}
}

@article{ogren2004cooperative,
  title={Cooperative control of mobile sensor networks: Adaptive gradient climbing in a distributed environment},
  author={Ogren, Petter and Fiorelli, Edward and Leonard, Naomi Ehrich},
  journal={IEEE Transactions on Automatic control},
  volume={49},
  number={8},
  pages={1292--1302},
  year={2004},
  publisher={IEEE}
}

@article{lynch2008decentralized,
  title={Decentralized environmental modeling by mobile sensor networks},
  author={Lynch, Kevin M and Schwartz, Ira B and Yang, Peng and Freeman, Randy A},
  journal={IEEE transactions on robotics},
  volume={24},
  number={3},
  pages={710--724},
  year={2008},
  publisher={IEEE}
}

@article{shriwastav2021dynamic,
  title={Dynamic Compressed Sensing of Unsteady Flows with a Mobile Robot},
  author={Shriwastav, Sachin and Snyder, Gregory and Song, Zhuoyuan},
  journal={arXiv preprint arXiv:2110.08658},
  year={2021}
}

@article{zhang2010cooperative,
  title={Cooperative filters and control for cooperative exploration},
  author={Zhang, Fumin and Leonard, Naomi Ehrich},
  journal={IEEE Transactions on Automatic Control},
  volume={55},
  number={3},
  pages={650--663},
  year={2010},
  publisher={IEEE}
}

@article{peng2014dynamic,
  title={Dynamic data driven application system for plume estimation using UAVs},
  author={Peng, Liqian and Lipinski, Doug and Mohseni, Kamran},
  journal={Journal of Intelligent \& Robotic Systems},
  volume={74},
  number={1},
  pages={421--436},
  year={2014},
  publisher={Springer}
}

@book{kutz2016dynamic,
  title={Dynamic mode decomposition: data-driven modeling of complex systems},
  author={Kutz, J Nathan and Brunton, Steven L and Brunton, Bingni W and Proctor, Joshua L},
  year={2016},
  publisher={SIAM}
}

@article{jovanovic2014sparsity,
  title={Sparsity-promoting dynamic mode decomposition},
  author={Jovanovi{\'c}, Mihailo R and Schmid, Peter J and Nichols, Joseph W},
  journal={Physics of Fluids},
  volume={26},
  number={2},
  pages={024103},
  year={2014},
  publisher={American Institute of Physics}
}

@article{askham2018variable,
  title={Variable projection methods for an optimized dynamic mode decomposition},
  author={Askham, Travis and Kutz, J Nathan},
  journal={SIAM Journal on Applied Dynamical Systems},
  volume={17},
  number={1},
  pages={380--416},
  year={2018},
  publisher={SIAM}
}

@article{Tu2014jcd,
	author = {J. H. Tu AND C. W. Rowley AND D. M. Luchtenburg AND S. L. Brunton AND J. N. Kutz},
	journal = {Journal of Computational Dynamics},
	number = {2},
	pages = {391--421},
	title = {On dynamic mode decomposition: theory and applications},
	volume = {1},
	year = {2014}}

@article{Brunton2022siamreview,
	author = {Brunton, Steven L and Budi{\v{s}}i{\'c}, Marko and Kaiser, Eurika and Kutz, J Nathan},
	journal = {SIAM Review},
	number = {2},
	pages = {229--340},
	title = {Modern {K}oopman Theory for Dynamical Systems},
	volume = {64},
	year = {2022}}

@article{Brunton2015jcd,
	author = {S. L. Brunton and J. L. Proctor and J. H. Tu and J. N. Kutz},
	journal = {Journal of Computational Dynamics},
	number = {2},
	pages = {165--191},
	title = {Compressed sensing and dynamic mode decomposition},
	volume = {2},
	year = {2015}}

@article{kramer2017sparse,
  title={Sparse sensing and DMD-based identification of flow regimes and bifurcations in complex flows},
  author={Kramer, Boris and Grover, Piyush and Boufounos, Petros and Nabi, Saleh and Benosman, Mouhacine},
  journal={SIAM Journal on Applied Dynamical Systems},
  volume={16},
  number={2},
  pages={1164--1196},
  year={2017},
  publisher={SIAM}
}

@article{kalman1960,
    author = {Kalman, R. E.},
    title = "{A New Approach to Linear Filtering and Prediction Problems}",
    journal = {Journal of Basic Engineering},
    volume = {82},
    number = {1},
    pages = {35-45},
    year = {1960},
    month = {03},
    issn = {0021-9223},
}

@article{gunnarson2021learning,
  title={Learning efficient navigation in vortical flow fields},
  author={Gunnarson, Peter and Mandralis, Ioannis and Novati, Guido and Koumoutsakos, Petros and Dabiri, John O},
  journal={Nature Communications},
  volume={12},
  number={1},
  pages={1--7},
  year={2021},
  publisher={Nature Publishing Group}
}

@article{krishna2022finite,
  title={Finite-horizon, energy-efficient trajectories in unsteady flows},
  author={Krishna, Kartik and Song, Zhuoyuan and Brunton, Steven L},
  journal={Proceedings of the Royal Society A},
  volume={478},
  number={2258},
  pages={20210255},
  year={2022},
  publisher={The Royal Society}
}

@article{biferale2019zermelo,
  title={Zermelo’s problem: optimal point-to-point navigation in 2D turbulent flows using reinforcement learning},
  author={Biferale, Luca and Bonaccorso, Fabio and Buzzicotti, Michele and Clark Di Leoni, Patricio and Gustavsson, Kristian},
  journal={Chaos: An Interdisciplinary Journal of Nonlinear Science},
  volume={29},
  number={10},
  pages={103138},
  year={2019},
  publisher={AIP Publishing LLC}
}

@inproceedings{buzzicotti2020optimal,
  title={Optimal Control of Point-to-Point Navigation in Turbulent Time Dependent Flows using Reinforcement Learning},
  author={Buzzicotti, Michele and Biferale, Luca and Bonaccorso, Fabio and Clark di Leoni, Patricio and Gustavsson, Kristian},
  booktitle={International Conference of the Italian Association for Artificial Intelligence},
  pages={223--234},
  year={2020},
  organization={Springer}
}

@article{madridano2021trajectory,
  title={Trajectory planning for multi-robot systems: Methods and applications},
  author={Madridano, {\'A}ngel and Al-Kaff, Abdulla and Mart{\'i}n, David and de la Escalera, Arturo},
  journal={Expert Systems with Applications},
  volume={173},
  pages={114660},
  year={2021},
  publisher={Elsevier}
}

@article{sargsyan2015nonlinear,
  title={Nonlinear model reduction for dynamical systems using sparse sensor locations from learned libraries},
  author={Sargsyan, Syuzanna and Brunton, Steven L and Kutz, J Nathan},
  journal={Physical Review E},
  volume={92},
  number={3},
  pages={033304},
  year={2015},
  publisher={APS}
}

@article{sargsyan2018online,
  title={Online interpolation point refinement for reduced-order models using a genetic algorithm},
  author={Sargsyan, Syuzanna and Brunton, Steven L and Kutz, J Nathan},
  journal={SIAM Journal on Scientific Computing},
  volume={40},
  number={1},
  pages={B283--B304},
  year={2018},
  publisher={SIAM}
}

@article{chaturantabut2010nonlinear,
  title={Nonlinear model reduction via discrete empirical interpolation},
  author={Chaturantabut, Saifon and Sorensen, Danny C},
  journal={SIAM Journal on Scientific Computing},
  volume={32},
  number={5},
  pages={2737--2764},
  year={2010},
  publisher={SIAM}
}

@book{Brunton2019book,
	author = {S. L. Brunton and J. N. Kutz},
	publisher = {Cambridge University Press},
	title = {Data-Driven Science and Engineering: Machine Learning, Dynamical Systems, and Control},
	year = {2019}}

@inproceedings{kalman1960general,
  title={On the general theory of control systems},
  author={Kalman, Rudolf E},
  booktitle={Proceedings First International Conference on Automatic Control, Moscow, USSR},
  pages={481--492},
  year={1960}
}

@article{reynolds2002improved,
  title={An improved in situ and satellite SST analysis for climate},
  author={Reynolds, Richard W and Rayner, Nick A and Smith, Thomas M and Stokes, Diane C and Wang, Wanqiu},
  journal={Journal of climate},
  volume={15},
  number={13},
  pages={1609--1625},
  year={2002},
  publisher={American Meteorological Society}
}

@dataset{noaasst,
  title={{NOAA Optimum Interpolation (OI) Sea Surface Temperature (SST) V2 data}},
  publisher={the NOAA PSL, Boulder, Colorado, USA},
  url={https://psl.noaa.gov}
}

\end{document}